\documentclass [12pt]{article}
\newfont{\bbb} {msbm10}
\newcommand{\Bbb}[1]{\mbox{\bbb#1}}
\newcommand{\R}{\Bbb{R}}
\newcommand{\bS}{\Bbb{S}}

\newcommand{\Z}{\Bbb{Z}}
\newcommand{\Q}{\Bbb{Q}}
\newcommand{\D}{\Bbb{D}}

\newcommand{\cT}{{\cal{T}}}
\newcommand{\sm}{\setminus}
\newcommand{\sbs}{\subset}
\newcommand{\ra}{\rightarrow}

\newcommand{\cP}{{\cal{P}}}

\newcommand{\cK}{{\cal{K}}}

\newcommand{\p}{\partial}

\newcommand{\cS}{{\cal{S}}}
\newcommand{\met}{{\cal{MET}}}
\newcommand{\mo}{{\cal{MET}}^{\, sec\, < \, 0}}

\newcommand{\To}{{\cal{T}}^{\, sec\, < \, 0}}

\newcommand{\cD}{{\cal{D}}}
\newcommand{\cC}{{\cal{C}}}

\newcommand{\cA}{{\cal{A}}}

\newcommand{\zp}{\Z_{(p)}}
\newcommand{\ozp}{\otimes\zp}
\usepackage[all]{xy}

\begin{document}

\title{ Teichm\"uller Spaces and Bundles with Negatively Curved Fibers}
\author{F. T. Farrell and P. Ontaneda\thanks{Both authors were
partially supported by NSF grants.}}
\date{}

\maketitle

\begin{abstract} 

In the first part of the paper we introduce the theory of bundles with negatively curved fibers. 
For a space $X$ there is a forgetful
map $F_X$ between bundle theories over $X$, which assigns to a bundle with negatively curved fibers over $X$ its subjacent smooth bundle.
Our Main result states that, for certain $k$-spheres $\,\bS^k$, the forgetful map $F_{\bS^k}$ is not one-to-one. This result follows from Theorem A, which proves that the quotient map \,{\scriptsize $ \mo(M)\ra\To(M)$}\, is not trivial at some homotopy levels,
provided the hyperbolic manifold $M$ satisfies certain conditions. Here \,{\scriptsize$\mo(M)$}\, is the space
of negatively curved metrics on $M$  and \,{\scriptsize $\To(M)=\mo(M)/DIFF_0(M)$}\, is, as defined in \cite{FO5}, the Teichm\"uller space of negatively curved metrics on $M$. Two remarks: {\bf (1)} the nontrivial elements in \,{\scriptsize $\pi_k\,\mo(M)$}\, constructed in \cite{FO6} have trivial
image  by the map induced by \,{\scriptsize $ \mo(M)\ra\To(M)$}\, {\bf (2)} the nonzero classes in \,{\scriptsize $\pi_k\,\To(M)$}\, constructed in \cite{FO5} are not in the image of the map induced by \,{\scriptsize $ \mo(M)\ra\To(M)$}\,; the nontrivial classes in \,{\scriptsize $\pi_k\,\To(M)$}\, given here, besides coming from \,{\scriptsize $\mo(M)$}\, and being
harder to construct, have a different nature and genesis: the former classes -given in \cite{FO5}- come from the existence of
exotic spheres, while the latter classes -given here- arise from the non triviality and structure of certain homotopy groups
of the space of pseudoisotopies of the circle $\bS^1$. The strength of the new techniques used here allowed us to prove also
a homology version of Theorem A, which is given in Theorem B.

\end{abstract}
\newpage

\noindent {\bf \Large  Section 0. Introduction.}\\

Let $M$ be a closed smooth manifold.
We will denote the group of all self diffeomorphisms of $M$, with the smooth topology, by $DIFF(M)$.
By a {\it smooth bundle over $X$, with fiber $M$}, we mean a locally trivial bundle for which the change of
coordinates between two local sections over, say, $U_\alpha, U_\beta\sbs X$ is given by a continuous map $U_\alpha\cap U_\beta\ra DIFF(M)$.
A {\it smooth bundle map} between two such bundles over $X$ is bundle map such that, when expressed in a local chart as
$U\times M\ra U\times M$, the induced map $U\ra DIFF(M)$ is continuous. In this case we say that the 
bundles are {\it smoothly equivalent}. 
Smooth bundles over a space $X$, with fiber $M$, modulo smooth equivalence, 
are classified by $\Big[ X, B\Big( DIFF(M) \Big)\Big]$, the set of homotopy classes of
(continuous) maps from $X$ to the classifying space $ B\Big( DIFF(M) \Big)$. \\

\noindent \Big[In what follows we will be considering everything {\it pointed}: $X$ comes with a base point $x_0$, the bundles come with
smooth identifications between the fibers over $x_0$ and $M$, and the bundle maps preserve these identifications. Also,
classifying maps are base point preserving maps.\Big]\\

If we assume that $X$ is simply connected, then
we obtain a reduction in the structural group of these bundles:
smooth bundles over a simply connected space $X$, with fiber $M$, modulo smooth equivalence,
are classified by $\Big[ X, B\Big( DIFF_0(M) \Big)\Big]$,
where $DIFF_0(M)$ is the space of all self diffeomorphisms of $M$ that are homotopic to
the identity $1_M$. In what follows we assume $X$ to be simply connected.
If we assume in addition that $M$ is aspherical with $\pi_1(M)$ centerless (e.g. admits a negatively curved metric) then
old results of Borel \cite{Borel}, Conner-Raymond \cite{CR} say that,  by pushing forward metrics,
$DIFF_0(M)$ acts freely on $\met(M)$, the space of smooth Riemannian metrics on $M$ (which we consider 
with the smooth topology). Moreover, Ebin's Slice Theorem \cite{Ebin} assures us that $DIFF_0(M)\ra \met(M)\ra\Big( \met(M)/DIFF_0(M)\Big)$
is a locally trivial bundle.
Hence, since $\met(M)$ is contractible, we can write $ B\Big( DIFF(M)_0 \Big)=\met(M)/DIFF_0$.
 In \cite{FO5} we called $\cT(M)=\met(M)\Big/ \Big(\R^+\times DIFF_0(M)\Big)$ the {\it Teichm\"uller Space of Riemannian Metrics on M},
where the $\R^+$ factor acts on $\met(M)$ by scalar multiplication. Since $\cT(M)$ is homotopy equivalent to $\met(M)/DIFF_0(M)$
we can also write $ B\Big( DIFF(M)_0 \Big)=\cT(M)$.
Therefore smooth bundles over a simply connected space $X$, with aspherical fiber $M$ and $\pi_1(M)$ centerless, 
modulo smooth equivalence, are classified by $\Big[ X, \cT(M)\Big]$.\\

Let $\cS$ be a complete collection of local sections of the bundle $\met(M)\ra\cT(M)$.
Using $\cS$ and a given map $f:X\ra\cT(M)$ we can explicitly construct a smooth bundle $E$ over $X$, with fiber $M$.
Yet, with these data we seem to get a little more: we get a Riemannian metric on each fiber $E_x$ of the bundle $E$.
This collection of Riemannian metrics does depend on $\cS$, but it is uniquely defined (i.e. independent of the choice
of $\cS$) up to smooth equivalence.\\

Of course, any bundle with fiber $M$ admits such a fiberwise collection of Riemannian metrics because $\met(M)$ is
contractible, so we seem to have gained nothing. On the other hand, in the presence a geometric condition we do get
a meaningful notion. We explain this next. \\

Denote by $\mo (M)$ the space
of all Riemannian metrics on $M$ with negative sectional
curvatures and by $\To(M)$ the image of $\mo(M)$ in $\cT(M)$ by the quotient map $\met(M)\ra\cT(M)$.
In \cite{FO5} we called $\To(M)$ the {\it Teichm\"uller Space of negatively Curved Riemannian Metrics on M}. 
If we are now given a map $X\ra\To(M)$, we get a smooth bundle $E$ with fiber $M$, and in addition, as mentioned before,
we get a collection of Riemannian metrics, one on each fiber $E_x$, $x\in X$. And, since now the target space is $\To(M)$, these Riemannian
metrics are all negatively curved. We call such a bundle a {\it bundle with negatively curved fibers}.
Still, to get a bona fide bundle theory we have to introduce the following concept. We say that two bundles $E_0$,
$E_1$ over $X$, with negatively curved fibers, are {\it negatively curved equivalent} if there is a bundle $E$ over
$X\times [0,1]$, with negatively curved fibers, such that $E|_{X\times\{ i\}}$ is smoothly equivalent to $E_i$, $i=0,1$,
via bundle maps that are isometries between fibers.
Then, bundles with negatively curved fibers over a (simply connected) space $X$, modulo negatively curved equivalence,
are classified by $\Big[ X, \To(M)\Big]$. And the inclusion map $F:\To(M)\hookrightarrow \cT(M)$ gives a relationship 
between the two bundle theories:\

$$ \Big[ X, \To(M)\Big]\stackrel{F_X}{\longrightarrow} \Big[ X, \cT(M)\Big]$$\

\noindent and the map $F_X$ is the ``forget the negatively curved structure" map.
The ``kernel" $\cK_X$ of this map between the two bundle theories is given by bundles over 
$X$, with negatively curved fibers, that are smoothly trivial.
Every bundle in $\cK_X$ can be represented by the choice of a negatively curved metric on each fiber of the
trivial bundle $X\times M$, that is, by a map $X\ra\mo(M)$. Note that this representation is not unique, because smoothly
equivalent representations give rise to the same bundle with negatively curved fibers. In any case, we have that $\cK_X$
is the image of $\Big[ X, \mo(M)\Big]$ by the map
$\Big[ X, \mo(M)\Big]\longrightarrow \Big[ X, \To(M)\Big]$,
induced by the quotient map $\mo(M)\ra\To(M)$.
Note that we can think of $\Big[ X, \mo(M)\Big]$ as a bundle theory: the ``bundles" here are choices of negatively curved
metrics, one for each fiber of the trivial bundle $X\times M$, modulo the following weak version of negatively curved
equivalence. Two ``bundles" $E_0$, $E_1$, here are equivalent if there is a ``bundle" $E$ over $X\times I$ such that
$E|_{X\times \{ i\}}=E_i$, $i=0,1$.
Summarizing, we get the following exact sequence of bundle theories:\vspace{.3in}

$(*)\hspace{.7in}\Big[ X, \mo(M)\Big]\stackrel{R_X}{\longrightarrow} \Big[ X, \To(M)\Big]\stackrel{F_X}{\longrightarrow} \Big[ X, \cT(M)\Big]$
\vspace{.3in}

\noindent where the map $R_X$ is the ``representation map", which, to each smoothly trivial bundle with negatively curved
fibers $E\in\cK_X$, assigns the set $R_X^{-1}(E)$ of representations of $E$ of the form $X\ra \mo(M)$.\\

It is natural to inquire about the characteristics of these maps. For instance, are they 
non constant? are they one-to-one? are they onto?
If, in (*), we specify $X=\bS^k$, $k>1$ (recall we are using basepoint preserving maps), we obtain $\pi_k(\mo(M))\ra\pi_k(\To(M))\ra\pi_k(\cT(M))$.
Some information about these maps between homotopy groups was given in \cite{FO5} and \cite{FO6}:\\

\begin{enumerate}
\item[1.] 
It was proved in \cite{FO6} that $\pi_2(\mo(M))$ is never trivial, provided $\mo(M)\neq\emptyset$ and $dim\, M>13$.
But the nonzero elements in  $\pi_2(\mo(M))$, constructed in \cite{FO6},  are mapped to zero by the map 
$\pi_2(\mo(M))\ra\pi_2(\To(M))$. Therefore the representation map $R_{\bS^2}$ in (*) is never one-to-one,
provided $\mo(M)\neq\emptyset$ and $dim\, M>13$.
\end{enumerate}

\noindent {\bf Remark.} It was also proved in \cite{FO6} (assuming $\mo(M)\neq\emptyset$) that $\pi_2(\mo(M))$ contains the infinite sum 
$(\Z_3)^\infty=\Big( \Z/3\Z \Big)^\infty$ as a subgroup, thus
$\pi_2(\mo(M))$ is not finitely generated. Moreover, it was proved that the same is true for $\pi_k(\mo(M))$,
for $k=2p-4$, $p>2$ prime (with $(\Z_p)^\infty$ instead of $(\Z_3)^\infty$) , provided $dim \, M$ is large (how large depending on $k$). 
Furthermore, $\pi_1(\mo(M))$ contains $(\Z_2)^\infty$, provided $dim\, M>11$.\\

\begin{enumerate}
\item[2.] The result of \cite{FO6} mentioned in the remark about $\pi_1(\mo(M))$ also proves that the
forget structure map $F_{\bS^2}$ is not onto. To see this just glue two copies of $\D^2\times M$ along $\bS^1$
using a nontrivial element in $\pi_1(\mo(M))$. Thus, there are (nontrivial) smooth bundles $E$ over $\bS^2$ which
do not admit a collection of negatively curved Riemannian metrics on the fibers of $E$. Using the remark, the same is true
for $\bS^k$, $k=2p-3$, $p>2$.

\item[3.] It was proved in \cite{FO5} that there are examples of closed hyperbolic manifolds
 for which $\pi_k(\To(M))$ is nonzero. Here $M$ depends on $k$ and always $k>0$. 
In \cite{FO5} no conclusion was reached on the case $k=0$ (i.e. about the connectedness of $\To(M)$).
Also, the images of these elements by the inclusion map $\To(M)\ra\cT(M)$ are not zero.
Hence the forget structure map $F_{\bS^k}$ is, in general,  not trivial. This means also that there are
bundles with negatively curved fibers that are not smoothly trivial, i.e. the representation map
$R_{\bS^k}$ is not onto in these cases.
\end{enumerate}

\noindent {\bf Remark.} In all the discussion above we can replace ``negatively curved metrics"
by ``$\epsilon$-pinched negatively curved
metrics".\\

Our main result here is the following:\\

\noindent {\bf Main Theorem.} {\it The forget structure map $F_{\bS^k}$ is, in general, not one-to-one, for
$k=2p-4$, $p $ prime.}\\

The Main Theorem follows from Theorems A, B and C, which actually prove more.
Theorems A and C together show that
for ``sufficiently large" closed hyperbolic $n$-manifolds the quotient map $\mo(M)\ra\To(M)$
is not trivial at the homotopy group level. That is,
$\pi_k(\mo(M))\ra \pi_k(\To(M))$ is  nonzero, 
provided a certain condition is satisfied by $k$ and $n$.
In particular this condition is satisfied for $k=0$ and $n>9$. Therefore we obtain as Corollary
that for sufficiently large closed hyperbolic $n$-manifolds, $n>9$, $\To(M)$
is disconnected. This solves the question left open in \cite{FO5} whether $\To(M)$ can ever be disconnected (see item 3 above). 
Also, the case $k=1$ is included here.\\

Theorem B and C together show that
for ``sufficiently large" closed hyperbolic $n$-manifolds the quotient map $\mo(M)\ra\To(M)$
is not trivial at the homology level. That is,
$H_k(\mo(M))\ra H_k(\To(M))$ is  nonzero, again
provided a certain condition is satisfied by $k$ and $n$. This is interesting because it gives characteristic classes
(mod a prime $p$) for the bundle theory. Finally, the case $k=1$ is also included here.
All results mentioned above also hold if we replace ``negative sectional curvature" by
``$\epsilon$-pinched to -1 sectional curvature".
To give more detailed statements of our results we need some notation.\\

Let $M$ be a closed hyperbolic manifold
and let $\gamma$ be an embedded closed geodesic in $M$.
We denote by    $\omega (\gamma)$ the width of its normal geodesic neighborhood. 
Given any $r>0$ and an embedded closed geodesic $\gamma$ in $M$
it is possible to find a finite sheeted cover $N$ of $M$
such that $\gamma$ lifts to a geodesic $\gamma_N$ in $N$ and $\omega(\gamma_N)>r$ (see \cite{FJ6}, Cor.3.3).\\

For a smooth closed manifold $N$ let $P(N)$ be the space of
topological radioisotopes of $N$, that is, the space of all homeomorphisms
$N\times I\ra N\times I$, $I=[0,1]$, that are the identity on $N\times \{ 0\}$.
We consider $P(N)$ with the compact open topology.
Also, $P^s(N)$ is the space of all {\it smooth} pseudoisotopies on $N$, with the smooth topology.
Let $TOP(N)$ be the group of self homeomorphisms of $N$, with the compact open topology. Note that
$DIFF(N)\sbs TOP(N)$.
We have the ``take top" map $\tau:P^s(N)\ra TOP(N)$, given by $\tau(f)=f|_{N\times\{ 1\}}:N\ra N$.\\

We will use the following notation. Let $X$ be a space. The geometric realization $|S(X)|$
of the singular simplicial set $S(X)$ of $X$ will be denoted by $X^\bullet$.
If $f:X\ra Y$ then we get the induced map $f^\bullet:X^\bullet\ra Y^\bullet$\\

Let $L\sbs TOP(\bS^1\times\bS^{n-2}))$
be the subgroup of all ``orthogonal" self  homeomorphisms of 
$\bS^1\times\bS^{n-2}$, that is $f:\bS^1\times\bS^{n-2}\ra\bS^1\times\bS^{n-2}$
belongs to $L$
if $f(z,u)=(e^{i\theta}z, A(z)u)$, for some $e^{i\theta}\in \bS^1$, and $A:\bS^1\ra SO(n-1)$.
We would like to use the quotient space $TOP(\bS^1\times\bS^{n-2}))/L$ but the quotient map $TOP(\bS^1\times\bS^{n-2})) 
\ra TOP(\bS^1\times\bS^{n-2}))/L\,\,\,$ is not a fibration. Instead we consider the {\it simplicial quotient} 
 \,\,\, $TOP(\bS^1\times\bS^{n-2}))//L:=\Big|\,S\Big(TOP(\bS^1\times\bS^{n-2})\Big)\Big/ S(L)\Big|$\,\,\, (see Section 1.3).
Now, define the map $\Upsilon _{n,\, k}:\pi_k (\, P^s(\bS^1\times\bS^{n-2})^\bullet\, )\ra \pi_k(\, TOP(\bS^1\times\bS^{n-2}) \, // \, L \,)$
as being the following composite:

$$\begin{array}{ccccc}
\pi_k (\, P^s(\bS^1\times\bS^{n-2})^\bullet\, ) & \stackrel{\pi_k(\tau^\bullet)}{\longrightarrow} & 
\pi_k(\, TOP(\bS^1\times\bS^{n-2})^\bullet \,) &
\longrightarrow & \pi_k(\, TOP(\bS^1\times\bS^{n-2}) \, //L\, ) 
\end{array}$$\\

\noindent where all homotopy groups have the corresponding identities as base points.
We will say that $\Upsilon _{n,\, k}$ is {\it strongly nonzero} if $\pi_k ( P^s(\bS^1\times\bS^{n-2})^\bullet)\cong
\pi_k ( P^s(\bS^1\times\bS^{n-2})\, )$
contains an infinite torsion subgroup $T$ such that $\Upsilon _{n,\, k}|_{\, T}$ is injective.\\

\noindent {\bf Theorem A.} {\it Let $k$ and $n> max\,\{ 3k+8,\, 2k+9 \}$ be such that $\Upsilon_{n,\, k}$ is strongly nonzero. 
Given $\ell>0$  there is a constant
$r=r(n,k, \ell)$ such that the following holds. If $M$ is a closed hyperbolic $n$-manifold 
that contains an embedded closed geodesic $\gamma$ with trivial normal bundle, length $\leq \ell$ and $\omega(\gamma)>r$,
then the map $\pi_k (\mo(M))\ra\pi_k(\To(M))$ is nonzero. In particular
$\pi_k(\To(M))$ is not trivial}.\\

In Theorem A all homotopy groups are based at (the class of) the given hyperbolic metric.
For $k=0$ the word ``nonzero" in the conclusion of Theorem A  should be read as ``not constant".
And the last sentence of Theorem A, when $k=0$, should be read as: ``$\To(M)$ is not connected".
Also, note that if $M$ is orientable, then condition ``$\gamma$ has trivial normal bundle" is redundant.\\

We have a homology version of Theorem A. Let ${\mit{h}}:\pi_k(\, TOP(\bS^1\times\bS^{n-2}) \, //L\, ) 
\ra H_k(\, TOP(\bS^1\times\bS^{n-2}) \, //\, L\, )$ denote the Hurewicz map.
As before, we will say that $h\,\Upsilon _{n,\, k}$ is {\it strongly nonzero} if $\pi_k ( P^s(\bS^1\times\bS^{n-2})\, )$
contains an infinite torsion subgroup $T$ such that $h\Upsilon _{n,\, k}|_{\, T}$ is injective.\\

\noindent {\bf Theorem B.} {\it Let $k>0$ and $n> max\,\{ 3k+8,\, 2k+9 \}$ be such that  \, h\,$\Upsilon_{n,\, k}$ is strongly nonzero. 
Given $\ell>0$  there is a constant
$r=r(n,k, \ell)$ such that the following holds. If $M$ is a closed hyperbolic $n$-manifold 
that contains an embedded closed geodesic $\gamma$ with trivial normal bundle, length $\leq \ell$ and $\omega(\gamma)>r$,
then the map $H_k (\mo(M))\ra H_k(\To(M))$ is nonzero. In particular
$H_k(\To(M))$ is not trivial}.\\

The statements of Theorems A  and B hold also for $\epsilon$-pinched negatively curved metrics:\\

\noindent {\bf Addendum to Theorems A and B.} {\it The statements of Theorems A and B remain true if we replace the 
decoration ``sec $<$ 0" on both $\mo(M)$ and $\To(M)$ by ``-1-$\epsilon<sec\leq$-1". But now
$r$ also depends on $\epsilon$, i.e.  $r=r(n,k,\ell, \epsilon )$.}\\

As mentioned before, 
the  condition ``$\gamma$ has a trivial normal bundle" in Theorems A and B can be obtained after taking, if
necessary, a two sheeted cover. The condition $\omega(\gamma)>r$ can also be obtained after taking a big enough
finite sheeted cover. To see this just take $r=r(n,k,\ell)$, $\ell=length (\gamma)$, and apply the result mentioned after the
the definition of $\omega(\gamma)$: {\it given any $r>0$ 
it is possible to find a finite sheeted cover $N$ of $M$
such that $\gamma$ lifts to a geodesic $\gamma_N$ in $N$ and $\omega(\gamma_N)>r$} (see \cite{FJ6}, Cor.3.3).
These facts imply the following results. \\

\noindent {\bf Corollary 1.} {\it Let $k>0$ and $n> max\,\{ 3k+8,\, 2k+9 \}$ be such that $h\Upsilon_{n,\, k}$ is strongly nonzero.
Then for every closed hyperbolic $n$-manifold $M$ there is a finite sheeted cover $N$ of $M$ such that
the maps $\pi_k (\mo(N))\ra\pi_k(\To(N))$, $H_k (\mo(M))\ra H_k(\To(M))$ are nonzero.}\\

And taking $k=0$ in Theorem A and Theorem C (see below) we have:\\


\noindent {\bf Corollary 2.} {\it 
Let $M$ be a closed hyperbolic $n$-manifold, $n>9$. Then $M$ admits a finite sheeted cover $N$ such that
$\To(N)$ is disconnected.}\\

\noindent {\bf Remark.}
The Addendum to Theorem A implies that the Corollaries remain true  if we replace the
decoration ``$sec\, <0$ " by ``-1-$\epsilon<sec\leq$-1". In this case $N$ depends not just on $n$ and $k$ but also on $\epsilon >0$.\\



Our next Theorem give cases for which $h\,\Upsilon_{n,\, k}$ is strongly not zero. We denote the
cyclic group of order $p$ by $\Z_p$, and $(\Z_p)^\infty$ is the (countably) infinite
sum of $\Z_p$'s. \\

\noindent {\bf Theorem C.} {\it Consider the map $h\Upsilon_{n,\, k}$. We have the following cases:}
\begin{enumerate}
\item[{}] {\bf k=0}\newline
{\it the group $\pi_0 ( P^s(\bS^1\times\bS^{n-2}))$ contains a subgroup $(\Z_2)^\infty$
and $h\Upsilon_{n,\, 0}$ restricted to this  $(\Z_2)^\infty$ is injective, provided $n\geq 10$.}

\item[{}] {\bf k=1}\newline
{\it the group $\pi_1 ( P^s(\bS^1\times\bS^{n-2}))$ contains a subgroup $(\Z_2)^\infty$
and $h\Upsilon_{n,\, 1}$ restricted to this  $(\Z_2)^\infty$ is injective, provided $n\geq 12$.}

\item[{}] {\bf k=2p-4, p$>$2 prime.}\newline
{\it the group $\pi_k ( P^s(\bS^1\times\bS^{n-2}))$ contains a subgroup $(\Z_p)^\infty$
and $h\Upsilon_{n,\, k}$ restricted to this  $(\Z_p)^\infty$ is injective, provided $n\geq 3k+8$.}

\end{enumerate}

\noindent {\bf Remark.} Of course if $h\Upsilon$ is strongly nonzero, then $\Upsilon$ is strongly nonzero, so
we have a statement similar to Theorem C for $\Upsilon$. In fact, this is Theorem D and it is used
to prove Theorem C (see Section 3.)

\newpage

\noindent {\bf \Large  Section 1. Preliminaries.}\vspace{.3in}


\noindent {\bf \large 1.1. Isotopies of metrics.}\\

For a Riemannian  manifold $Q$, with metric $g$,  and a submanifold $P$ we denote by $\perp_pP$ the orthogonal complement of
$T_pP$ in $T_pQ$, with respect to the metric $g$. As usual the exponential map $TQ\ra Q$ is denoted by $exp$ and, to avoid complicating
our notation, the normal (to $P$) exponential map will also be denoted: $exp:\perp_pP\ra Q$.
If we need to show the dependence of these objects on $g$ we shall write
$\perp^gP$ and $exp^{\, g}$. Recall that if $E$ is a subbundle of $TQ|_P$ such that $E\oplus TP=TQ|_P$ then
$exp|_E:E\ra Q$ (here $exp:TQ|_P\ra Q$) is a diffeomorphism near $P$, is the identity on $P$, and the derivative at any point
 $p\in P$ is the identity
(after the obvious identification of $TE|_P$ with $TQ|_P$).
We will need the following result.\\

\noindent {\bf Proposition 1.1.1.} {\it Let $M$ be smooth $n$-manifold without boundary,
and $g_0$, $g_1$ two Riemannian metrics on $M$.
Let also $P$ be a closed smooth $k$-submanifold of $M$, $2k+3\leq n$ with trivial normal bundle 
and $\eta : P\ra M$ a smooth embedding homotopic to the
inclusion $\iota :P\hookrightarrow M$. Then there is a smooth isotopy $h_t:M\ra M$, $0\leq t\leq 1$, $h_0=1_M$,
such that (write $h=h_1$, and $g=h^*g_1$):}
\begin{enumerate}
\item[{1.}] {\it $\eta=h\iota$.}

\item[{2.}] {\it $\perp^g_pP= \perp^{g_0}_pP$, for all $p\in P$. }

\item[{3.}] {\it $g(u,v)=g_0(u,v)$, for all $u,v\in \perp^g_pP$ .}

\item[{4.}] {\it  There is $\epsilon>0$ such that $exp^{\, g}_p(v)=exp^{\, g_0}_p(v)$, for all
$p\in P$, $v\in \perp^g_pP$ with $g(v,v)\leq\epsilon$.}
\end{enumerate}

\noindent{\bf Proof.} If $\iota=\eta$ and $g_0=g_1$  we are done. If not let $H:P\times [0,1]\ra M$ be a homotopy between $\iota$ and $\eta$.
Since $2k+3\leq n$ we can assume that $H $ is an embedding. Hence $\iota $ and $\eta$ are
isotopic. This isotopy can be extended in the usual way (using vector fields) to an ambient
isotopy of $M$. In this way obtain an isotopy that satisfies (1). We will construct other
isotopies to obtain (2)-(4) (these remaining isotopies will fix $P$). Hence we can assume that
$\iota=\eta$.\\

Let $V=\{ v_1, ...,v_l\}$, $l=n-k$, be an orthonormal framing of the bundle $\perp^{\,g_0}P\sbs TM|_P$.
Let also $V'=\{ v'_1, ...,v'_l\}$ be the projection of $V$ in $\perp^{\,g_1}P$, that is 
$v'_i=u_i+v_i$, $u_i\in TP$ and $v'_i\in \perp^{\,g_1}P$.
Since $\perp^{\, g_1}P\oplus TP=TM$ we have that $V'$ is a framing of $\perp^{\,g_1}P$.
Denote by $\Phi_t:\perp^{g_0}P\ra TM|_P$ the bundle map given by $\Phi_t(v_i)=tu_i+v_i$, and
let $E_t$ be the subbundle of $TM|_P$ generated by the $tu_i+v_i$, that is $E_t=\Phi_t(\perp^{g_0}P)$.
Then $E_t\oplus TP=TM|_P$. Let $exp^t :E_t\ra M$ be the restriction of $exp:TM|_P\ra M$.
Then $H_t=exp^t\circ \Phi_t\circ(exp^0)^{-1}$ is an isotopy defined on a neighborhood of $P$, starting at the identity.
And, since the derivative of $exp^t$ at a $p\in P$ is the identity we have that the derivative of $H_1$ 
at a $p\in P$ sends $\perp^{g_0}P$ to $\perp^{g_1}P$. Extend $H_t$ to the whole $M$. It is not difficult
to show that $(H_1)^*g_1$ satisfies item (2) (with $(H_1)^*g_1$ instead of $g$). Hence we
can suppose now that $\perp^{g_1}_pP= \perp^{g_0}_pP$, for all $p\in P$.\\

We now further change $g_1$ by an isotopy so as to obtain (3). Note that if $V$ is also
orthonormal with respect to $g_1$ we are done. If not let $V_t$ be a path of framings of
$\perp^{g_1}P= \perp^{g_0}P$ with $V_0=V$ and $V_1$ orthonormal with respect to $g_1$
(for this just apply the canonical Gram-Schmidt orthonormalization process). Let 
$\Phi'_t:\perp^{g_0}P\ra \perp^{g_0}P$ be the bundle map that sends $V$ to $V_t$. Let
$H'_t=exp\circ \Phi'_t\circ(exp)^{-1}$ (here $exp:\perp^{g_0}P\ra M$) is an isotopy defined on a neighborhood of $P$, starting at the identity.
And the derivative of $H_1$ 
at a $p\in P$ sends $V$ to $V_1$. Extend $H_t$ to the whole $M$. It is not difficult
to show that $(H'_1)^*g_1$ satisfies item (3) (with $(H_1)^*g_1$ instead of $g$). Hence we
can suppose now that  $g_1(u,v)=g_0(u,v)$, for all $u,v\in \perp^g_pP$.\\

Now, note that the map 
$exp^{ \,\, g_1}\circ (exp^{\,\, g_0})^{-1}$ is a diffeomorphism defined on a neighborhood of $P$, and
its derivative at a $p\in P$ is the identity. Hence a fiber version of Alexander's trick (see Appendix)
gives an isotopy $H''_t$ that deforms
$exp^{\,\, g_1}\circ (exp^{\,\, g_0})^{-1}$ to the identity (near $P$). Extending this isotopy to the whole
$M$, we have that  $(H''_1)^*g_1$ satisfies item (4). This proves the Proposition.\\

\noindent {\bf Remarks.} 

\noindent {\bf 1.} From the proof of the Proposition we see that the isotopy $h_t$ in the
statement of the Proposition can be chosen to have
support in a small (as small as we want) neighborhood of the image of the embedding $H$ mentioned in the first
paragraph of the proof.

\noindent {\bf 2.} 
It can be checked from the proof above that a (local) parametrized version of Prop. 1.1.1 also holds:
fix $g_0$ and suppose that we are given a $C^2$-neighborhood  $U$ of some $g_1$ in $\met (M)$,
and a continuous map $g'\mapsto \eta_{g'}\in Emb^{\infty}(P, M)$, $g'\in U$. Then the proof of Prop. 1.1.1 above gives us
a method to construct a map $h_{g'}$
such that the map $g'\mapsto h_{g'}$ (and hence the map
$g'\mapsto h^*_{g'}(g')$) is well defined and continuous on some $C^2$-open $W\sbs U$, $g_1\in W$. 
(Here $\epsilon$ will depend on  $W$.) 
Moreover, if $g'$ already satisfies (1)-(4) of the Proposition then 
this map leaves $g'$ invariant, that is, $g'=h^*_{g'}(g')$.
\vspace{.5in}

\noindent {\bf \large 1.2. The map $\Omega_\gamma$}\\

We will need the following construction, which have some similarities to the one given in
\cite{FO}. \\

Write $\bS^1(\ell )=\{ (x,y)\in\R^2\, ,\,\, x^2+y^2=(\ell/2\pi)^2\}$. Let $M$ be a hyperbolic
$n$-manifold, with metric $g_0$. Let $\gamma :\bS^1(\ell)\ra M$ be an embedded closed geodesic of length $\ell$. 
Sometimes we will denote the image of $\gamma $ just by $\gamma$.
We assume that the normal bundle of $\gamma$ is orientable, hence trivial. 
Let $r>0$ be such that $6r$ is less than the width of the normal geodesic tubular
neighborhood of $\gamma$ and denote by $U$ the normal geodesic tubular neighborhood of $\gamma$ of width 6$r$.
Write $V=U\sm \gamma$. 
Using the exponential map of geodesics orthogonal
to $\gamma$ and parallel translation along $\gamma$ we get that $V$ (with the given hyperbolic metric $\rho_0$)
is isometric to the quotient of $\R\times\bS^{n-2}\times (0,6r]$, equipped with the doubly warped Riemannian metric:

$$\rho'(s,u,t)=\cosh^{2}(t)ds^2+\sinh^{2}(t) \sigma_{\bS^{n-2}}(u)+dt^{2},$$

\noindent  by the action of an isometry $A:\R\times\bS^{n-2}\times (0,6r]\ra\R\times\bS^{n-2}\times (0,6r]$
of the form $(s,u,t)\mapsto (s+\ell, Tu, t)$ for some $T\in SO(n-1)$.
Here $\sigma_{\bS^k}$ is the canonical round Riemannian metric on the $k$-sphere $\bS^k$.\\

\noindent {\bf Remark.} Note that $V$ is diffeomorphic to $\bS^1\times\bS^{n-2}\times (0,6r]$. For a
compactness argument that will be used later we need some canonical ways of identifying
$V$ with $\bS^1\times\bS^{n-2}\times (0,6r]$. We do this by choosing certain $A$-invariant 
trivializations of the bundle  $\R\times\R^{n-1}\ra \R$. To do this consider $SO(n-1)$ with its bi-invariant
metric and let $B_1,...,B_k$ be closed geodesic balls that cover $SO(n-1)$ and such that:
for each $S\in B_i$ there is a path $\alpha_{S,\, i}\, (u)$, $u\in [0,1]$, starting at the identity and ending
in $S$. Also, for each $i$ fixed, $\alpha_{S,\, i}$ varies continuously (in the smooth topology) with $S$.
We require also that $\alpha_{S,\, i}$ be constant near 0 and 1.
(For instance, $\alpha_{S,\, i}$ could be a properly rescaled geodesic.)
Then our canonical $A$-invariant trivializations are constructed in the following way.
For $T\in B_i$ define  $v^i_j(t)=\alpha_{T,\, i}\, (t/\ell ).\, e_j$, $j=1,...,n-1$, for $t\in [0,\ell]$
(where the $e_j$'s form the canonical base of $\R^{n-1}$). Note that $A(0,v^i_j(0))=A(0, e_j)=
(\ell, Te_j)=(\ell, v^i_j(\ell))$,
thus we can extend the $v^i_j(t)$ periodically to all $t\in \R$. Therefore, for each $i$ such that $T\in B_i$,  $\Big\{ v^i_1,...,v^i_{n-1}\Big\}$
is a $A$-invariant trivialization. 
Hence for each $T$ we get finitely many ``canonical"
trivializations, one for each $i$ such that $T\in B_i$. \\

\noindent{\it Caveat:} We are giving two identifications of the universal cover of $V$ with 
$\R\times\bS^{n-2}\times (0,6r]$: (1) $V$ is the quotient of $\R\times\bS^{n-2}\times (0,6r]$  by the action of the isometry $A$,
and (2) using the canonical trivializations mentioned in this remark. These two identifications do not necessarily
coincide.\\

Let $\delta:[0,1]\ra [0,1]$ be a smooth map such that $\delta(0)=0$, $\delta(1)=1$ and which is constant
near 0 and 1.
Define the metric $\rho''$ on $\R\times\bS^{n-2}\times (0,6r]$ in the following way:\\
\begin{enumerate}
\item[{$\bullet$}] $\rho''=\rho'$  outside $\R\times \bS^{n-2}\times [2r,5r]$.

\item[{$\bullet$}] On 
$\R\times \bS^{n-2}\times [2r,3r]$ we have:
$$\rho''(s,u,t)=\cosh^{2}(t) ds^2+\frac{1}{4} \left[ e^t+\left( 2\,\delta(\frac{t-2r}{r})-1\right)
e^{-t}\right]^{2} \sigma_{\bS^{n-2}}(u)+dt^{2}$$
\noindent Note that for $t\leq 3r$, and near $3r$, we have $\rho''(s,u,t)=\cosh^{2}(t)[ ds^2+ \sigma_{\bS^{n-2}}(u)]+dt^{2}$.
That is, $\rho''$ is a {\it simply} warped metric in this case.

\item[{$\bullet$}] On 
$\R\times \bS^{n-2}\times [3r,4r]$ define:
$\rho''(s,u,t)=\cosh^{2}(t)[ ds^2+ \sigma_{\bS^{n-2}}(u)]+dt^{2}$.

\item[{$\bullet$}] On 
$\R\times \bS^{n-2}\times [4r,5r]$ define:
$$\rho''(s,u,t)=\cosh^{2}(t) ds^2+\frac{1}{4} \left[ e^t+\left( 1-2\,\delta(\frac{t-4r}{r})\right)
e^{-t}\right]^{2} \sigma_{\bS^{n-2}}(u)+dt^{2}$$

\end{enumerate}

Note that $\rho''$ is also invariant by $A$, hence induces a Riemannian metric $\rho_1$ on
the quotient $V$. \\

\noindent {\bf Lemma 1.2.1.} {\it Given $\epsilon>0$, we have that all sectional curvatures of
$\rho''$ and $\rho_1$ lie in the interval $(-1-\epsilon, -1+\epsilon)$, provided $r$
is sufficiently large (how large depending solely on $\epsilon$ and $n$).}\\

\noindent For a proof see Lemma 1.2.1 in \cite{FO}.\\

\noindent {\bf Remark.} In the next Section we will need a canonical way of deforming $\rho''$ to $\rho'$.
To do this we assume in addition that $\delta(1-v)=1-\delta(v)$, for $v\in [0,1]$, and define
$\rho''_v$ in the following way. Define $\rho''_v=\rho''$ on $\R\times\bS^{n-2}\times\Big([2r,2r+vr]\cup [5r-vr,5r] \Big)$
and $\rho''_v(s,u,t)=\cosh^{2}(t) ds^2+\frac{1}{4} \Big[ e^t+( 2\,\delta(v)-1 )e^{-t}\Big]^{2} 
\sigma_{\bS^{n-2}}(u)+dt^{2}$ on $\R\times\bS^{n-2}\times [2r+vr, 5r-vr]$.
(And also let $\rho''_v=\rho''=\rho'$ outside $\R\times\bS^{n-2}\times [2r,5r]$.)
Then $\rho''_0=\rho'$ and $\rho''_1=\rho''$. Also $\rho''_v$ is invariant by $A$, hence induces a deformation $(\rho_1)_v$ 
from $\rho_0$ to $\rho_1$. Furthermore, Lemma 1.2.1 holds for all $\rho''_v$ and $(\rho_1)_v$, provided $r$
is sufficiently large.\\

We now define  a map $\Omega_\gamma : P^s(\bS^{1}\times\bS^{n-2})\ra \met (M)$.
Identify $V$ with $\bS^1\times \bS^{n-2}\times (0,6r]$
(via one of the finitely many canonical ways mentioned in the remark above). Note that for
the metric $\rho_1$ on $V=\bS^1\times \bS^{n-2}\times (0,6r]$ we have that: {\bf (1)}
$\rho_1=\rho_0$  outside $\bS^1\times \bS^{n-2}\times [2r,5r]$ and {\bf (2)} on
$\bS^1\times \bS^{n-2}\times [3r,4r]$ we have
$\rho_1(z,u,t)=\cosh^{2}(t)\sigma(z,u)+dt^{2}$ for some metric $\sigma$ on $\bS^1\times\bS^{n-2}$.
For $\varphi\in P^s(\bS^1\times\bS^{n-2} )$ 
write $\mu=\tau\varphi\times1_{[0,1]}\in DIFF(\bS^1\times\bS^{n-2}\times [0,1])$ and
let $\omega=\Omega_\gamma (\varphi )$  be the metric 
defined as follows:
\begin{enumerate}
\item[{\bf a.}] $\omega$ is the given hyperbolic metric $\rho_0$ (or $\rho_1$) outside $\bS^1\times \bS^{n-2}\times [r,5r]\sbs M$.
And $\omega=\rho_1$ on $\bS^1\times \bS^{n-2}\times [4r,5r]\sbs M$

\item[{\bf b.}] $\omega= [\lambda^{-1}\varphi\lambda]^*\rho_0=[\lambda^{-1}\varphi\lambda]^*\rho_1$,
where $\lambda(z,u,t)=(z,u,\frac{t-r}{r})$, for $(z,u,t)\in\bS^1\times \bS^{n-2}\times [r,2r]$.

\item[{\bf c.}] $\omega=[\lambda^{-1}\mu\lambda]^*\rho_1$,
where $\lambda(z,u,t)=(z,u,\frac{t-2r}{r})$, for $(z,u,t)\in\bS^1\times \bS^{n-2}\times [2r,3r]$.
Note that for $t\leq 3r$, and $t$ near $3r$, 
$\omega (z,u,t)=\cosh^{2}(t)[(\tau\varphi)^*\sigma(z,u)]+dt^{2}$ 

\item[{\bf d.}] On 
$\bS^1\times \bS^{n-2}\times [3r,4r]$ define:
$$\omega(z,u,t)=\cosh^{2}(t) \sigma_{t}(u)+dt^{2}$$
\noindent where $\sigma_t=\left(1-\delta(\frac{t-3r}{r})\right)(\tau\varphi)^*\sigma (z,u)+\delta(\frac{t-3r}{r})\sigma(z,u)$.
\end{enumerate}

\noindent{\bf Remark.} 
We are assuming that all pseudoisotopies are products near 0 and 1. Hence items
(a)-(d) give a well defined Riemannian metric $\omega$ on $M$.\\

\noindent {\bf Lemma 1.2.2.} {\it Given $\epsilon >0$, $\ell>0$ and a compact set $K\sbs P^s(\bS^1\times \bS^{n-2})$, there is $r>0$ such that
the sectional curvatures of $\, \Omega_\gamma (\varphi)$ lie in $(-1-\epsilon , -1+\epsilon) $, for all $\varphi\in K$, provided
$\gamma $ has length $\ell$.}\\

\noindent {\bf Proof.} Outside $\bS^1\times \bS^{n-2}\times [3r,4r]$, $\omega=\Omega_\gamma (\varphi)$ coincides with
$\rho_1$ or a pullback of it. Hence Lemma 1.2.1 implies in this case all
sectional curvatures of $\, \Omega_\gamma (\varphi)$ lie in $(-1-\epsilon , -1+\epsilon) $.
For $\bS^1\times \bS^{n-2}\times [3r,4r]$ we apply Lemma 2.2 of \cite{O}. 
How large $r$ needs to be in this case
depends only on $\epsilon$, $(\tau\varphi)^*\sigma$
and their derivatives up to order 2. Since $K$ is compact $(\tau\varphi)^*\sigma$
and their derivatives up to order 2 are bounded provided that all possible
metrics $\sigma$ on $\bS^1\times \bS^{n-2}$ and their derivatives up to order 2 are bounded.
Recall that these metrics $\sigma$ are obtained in the following way.   For $T\in SO(n-1)$ let $\alpha_{T,\, i}: [0,1]\ra SO(n-1)$,
with $T\in B_i$ be as in the remark above. Let $f:[0,1]\times\bS^{n-2}\ra [0,\ell ]\times\bS^{n-2}$
defined by $f(t,u)=(\ell t, \alpha_{T,\, i}(t).u)$. Let $\sigma'=f^*(ds^2+\sigma_{\bS^{n-2}})$, where
$ds^2+\sigma_{\bS^{n-2}}$ is the canonical product metric on $[0,\ell ]\times\bS^{n-2}$. Gluing $\{ 0\}\times\bS^{n-2}$
to $\{ 1\}\times\bS^{n-2}$ the metric $\sigma'$ gives a metric $\sigma$ on $\bS^1\times\bS^{n-2}$.
Since $\ell$ is fixed and the set of all $\alpha_{T,\, i}$ is compact in $DIFF([0,1],SO(n-1))$ (with the smooth topology),
it follows that all possible
metrics $\sigma$ on $\bS^1\times \bS^{n-2}$ and their derivatives up to order 2 are bounded.
This proves the Lemma.\\

\noindent {\bf Remark.} A subtle point here. In Lemma 1.2.2 the number $r$ depends only on $\epsilon$, $\ell$, $K$ and
the dimension $n$ of the manifold, but not on the particular manifold $M$. The independence from $M$ stems from the canonical
identifications mentioned in the remark at the beginning of this Section. 
\vspace{.8in}

\noindent {\bf 1.3. The Simplicial Quotient.}\\

Let $X$ be a space and $S(X)$ be its singular simplicial set. Recall that the $q$-simplices of $S(X)$
are the singular $q$-simplices on $X$, i.e maps $\Delta^q\ra X$. Write $X^\bullet=|\, S(X)|$
where the bars denote ``geometric realization". There is a canonical map $h_X:X^\bullet\ra X$
which is a weak homotopy equivalence. If $f:X\ra Y$ is a map then the simplicial map 
$S(f) :S(X)\ra S(Y)$ defines a map $f^\bullet :X^\bullet\ra Y^\bullet$ and clearly

$$\begin{array}{ccc}
X^\bullet&\stackrel{f^\bullet}{\ra}& Y^\bullet\\
\downarrow&&\downarrow\\
X&\stackrel{f}{\ra} & Y\end{array}$$

\noindent commutes.\\

Let $G$ be a topological group
acting freely on $X$. Then $S(G)$ is a simplicial group acting simplicially on $S(X)$
and we get a simplicial set $S(X)/S(G)$. We define the {\it simplicial quotient} as $X//G=|S(X)/S(G)|$. 
The map $S(X)\ra S(X)/S(G)$ defines a map $X^\bullet\ra X//G$.  
We will use the following facts.\\

\noindent {\bf 1.} We have that $X^\bullet\ra X//G$ is a fibration with fiber $L^\bullet$.\\

\noindent {\bf 2.} Let $G$ and $H$ act freely on $X$ and $Y$, respectively. Let $f:X\ra Y$ be continuous
and $D:G\ra H$ be a homomorphism (or anti homomorphism) of topological groups, and assume that $f$ is $D$-equivariant, that is,
$f(gx)=D(g)f(x)$, for all $x\in X$ and $g\in G$. Then $f$ defines a map $F:X//G\ra Y//H$ and

$$\begin{array}{ccc}
X^\bullet&\stackrel{f^\bullet}{\ra}& Y^\bullet\\
\downarrow&&\downarrow\\
X//G&\stackrel{F}{\ra} & Y//H\end{array}$$

\noindent is commutative.\\

\noindent {\bf 3.} 
If $q:X\ra X/G$ is a (locally trivial) fiber bundle then the simplicial map $S(X)\ra S(X/G)$ is onto. Furthermore,
two singular simplices in $X$ have the same image in $S(X/G)$ iff they differ by an element in
$S(G)$. Hence the simplicial map $S(X)/S(G)\ra S(X/G)$ is a bijection. It follows that
$X//G\ra (X/G)^\bullet$ is a homeomorphism.
\vspace{.4in}

\noindent {\bf 1.4. The space $\met(Q,g)$.}\\

We have considered the space of Riemannian metrics of a closed manifold. We now mention some facts and give a
few definitions related to the non compact case. Let $Q$ be a complete Riemannian manifold, with metric $g$.
We consider $\met(Q)$ to be the set of complete Riemannian metrics on $Q$ with the smooth topology, which is the union, 
for all $k$, of the topologies of
$C^k$-convergence on compact sets. Similar topology is given to $DIFF(Q)$.\\

Let $f:(X_1,d_1)\ra (X_2,d_2)$ be surjective map between metric spaces. Recall that, in this particular case,
$f$ is a $(\lambda ,\delta)$-quasi-isometry
if $\frac{1}{\lambda}\, d_1(x,y)-\delta\leq d_2(f(x),f(y))\leq \lambda\, d_1(x,y)+\delta$, for all
$x,y\in X_1$. 
Let $g', g''\in\met(Q)$. We say that $g'$ and $g''$ are $(\lambda ,\delta)$-quasi-isometric if the identity 
$(Q, d_{g'})\ra(Q,d_{g''})$ is a  $(\lambda ,\delta)$-quasi-isometry, where $ d_{g'}$, $d_{g''}$ are the intrinsic 
metrics induced by the Riemannian metrics
$g'$ and $g''$, respectively. A useful way to prove that two metrics are quasi-isometric is the following: \\

\noindent {\bf 1.4.1. Lemma.} {\it Let $g,\, g'$ be two complete Riemannian metrics on the manifold $Q$.
Suppose there are constants $a,b>0$ such that
$a^2\leq g'(w,w)\leq b^2$ for every $w\in TQ$ with $g(w,w)=1$. Then $g$ and $g'$ are
$(\lambda,0)$-quasi-isometric, where $\lambda=max\{ \frac{1}{a},\, b\}$.}\vspace{.1in}

The proof is straightforward (see Lemma 2.1 of \cite{FO6}).\\

Here is a variation of the space $\met(Q)$. We define $\met(Q,g)$ to be the set of complete Riemannian metrics on $Q$ that are
quasi-isometric to $g$. We give $\met(Q,g)$ the smooth quasi-isometry topology: basic neighborhoods of a $g'\in\met(Q,g)$ are intersections of
open neighborhoods (in $\met(Q)$) of $g'$ with the quasi-geodesic balls $$B_{\lambda , \delta}(g')=\{h\in\met(Q,g)\, /\, h\,\,{\mbox{is }} (\lambda', \delta')-{\mbox{quasi-isometric to}}\,\, g',\,\, \lambda'<\lambda,\,\delta'<\delta\}$$
Then the inclusion $\met(Q,g)\hookrightarrow \met(Q)$ is continuous, but the topology of $\met(Q,g)$ is strictly finer than the one induced by $\met(Q)$.\\

Let $M$ be closed and let $p:Q\ra M$ be a covering map. Let $g\in\met(Q)$ be such that $g$ is quasi-isometric to a mertric of the form
$p^*(g')$ for some (hence all) $g'\in\met(M)$.
Then the map
$$\met(M)\stackrel{{\mbox{lift}}}{\longrightarrow}\met(Q,g)$$
given by $g'\mapsto p^*(g')$ is well defined and clearly continuous. For instance we can take $g=p^*(g')$, $g'\in\met(M)$.
Note that the topology of $\met(Q,p^*(g'))$ is independent of the choice of $g'$.\\
\\

Let $DIFF(Q,g)$ be the subset of $DIFF(Q)$ of all self-diffeomorphisms $\phi$ that are at bounded $g$-distance from the identity $1_Q$;
that is $d_g(\phi, 1_Q)=sup\{ d_g(\phi (x), x)\, /\, x\in Q\}$ is finite.
We give $DIFF(Q,g)$ the smooth quasi-isometry topology: the open basic sets of the identity $1_Q$ are intersections of open neighborhoods of
$1_Q$ in $DIFF(Q)$ with the sets $\{ \phi\, /\,  d_g(\phi,1_Q)<\epsilon  \}$.
We define $DIFF_0(Q,g)$ to be the subspace of $DIFF(Q,g)$ of self-diffeomorphisms of $Q$ that are $g$-boundedly homotopic to $1_Q$. \\

We have that 
the action of $DIFF(Q,g)$ on $\met(Q,g)$ is continuous.\\ 

Of course if $Q$ is closed then $\met(Q)$ coincides with $\met(Q,g)$, for any $g$.\\

Finally, define $\To(Q,g)=\mo(M,g)/\R^+\times DIFF_0(Q,g)$.
\newpage

\noindent {\bf 1.5. The Space at Infinity.}\\

Let $H$ be a complete, simply connected manifold of nonpositive curvature, that is, $H$ is a Hadamard manifold.
Recall that the space at infinity $\p_\infty H$ is defined as the quotient of the set of geodesic rays by the relation:
``finite Hausdorff distance'' (see, for instance, \cite{BH}).
In this definition ``geodesic rays'' can be replaced by ``quasi-geodesic rays'',
provided $H$ has sectional curvatures $\leq c <0$. The compactification
${\overline{H}}=H\cup\p_\infty H$ is given the ``cone topology". We mention three useful facts:

\begin{enumerate}
\item[1.] The definition of the cone topology implies the following. If $\beta$ is a geodesic ray between $p\in H$ and
$q\in\p_\infty H$, and $V$ is a neighborhood of $q$ in $\overline{H}$, then there is $T>0$ such that the following holds.
For any other geodesic ray $\beta'$ such that the distance between $\beta(t)$ and  $\beta'(t)$ is
$\leq 1$, for $t\in[0,T]$,  we have that $\beta'(t)\in V$, for all $t\in [T,\infty ]$.

\item[2.] If we assume that $H$ has sectional curvatures $\leq c <0$, then we get the following
quasi-geodesic version of item 1. Let $\beta$ be a geodesic ray between $p\in H$ and 
$q\in\p_\infty H$, $V$ is a neighborhood of $q$ in $\overline{H}$ and $\lambda>0$, $\delta\geq 0$.
Then there is $T>0$ such that 
for any $(\lambda,\delta)$-quasi-geodesic ray $\beta'$ for which the distance between $\beta(t)$ and  $\beta'(t)$ is
$\leq 1$, for $t\in[0,T]$, we have that $\beta'(t)\in V$, for all $t\in [T,\infty ]$.

\item[3.] 
If $g_1$ and $g_2$ are  quasi-isometric complete Riemannian metrics on 
the simply connected manifold $H$, with sectional curvatures $\leq c <0$, then the space at infinity
and the compactification of $H$ are the same (as topological spaces) if taken with respect to $g_1$ or $g_2$.

\end{enumerate}

We can generalize most of the concepts mentioned above to the following non simply connected case (see Section 2 of \cite{FO6}). 
Let $Q$ be a complete Riemannian manifold (with metric $g$) with sectional curvatures $\leq c <0$.
Let also $S$ be a closed totally geodesic submanifold of $Q$ such that
$\pi_1(S)\ra\pi_1(Q)$ is an isomorphism. Then $Q$ is diffeomorphic to the total space of the normal of $S$ in $Q$, via
the normal (to $S$) exponential map. A geodesic ray (i.e. a local isometry $[0,\infty)\ra Q$) 
either diverges from $S$ or stays at bounded distance from
$S$. Then the space at infinity of $Q$ can be defined as before:
the space at infinity $\p_\infty Q$ is defined as the quotient of the set of geodesic rays that diverge from $S$, by the relation:
``finite Hausdorff distance''. In this definition we can replace ``geodesic rays that diverge from $S$" by `` quasi-geodesic rays".
The compactification is
${\overline{Q}}=Q\cup\p_\infty Q$ and is given the ``quotient cone topology". Then, in this context, we also get 
(almost) exact versions 1',2' and 3'
of items 1, 2, and 3 above.
We will use only use 2' and 3'. Here they are (see Section 2 of \cite{FO6}).

\begin{enumerate}
\item[2'.] Let $\beta$ be a geodesic ray between $p\in Q$ and 
$q\in\p_\infty Q$, $V$ a neighborhood of $q$ in $\overline{Q}$ and $\lambda>0$, $\delta\geq 0$.
Then there is $T>0$ such that 
for any $(\lambda,\delta)$-quasi-geodesic ray $\beta'$ for which the distance between $\beta(t)$ and  $\beta'(t)$ is
$\leq 1$, for $t\in[0,T]$,
we have that $\beta'(t)\in V$, for all $t\in [T,\infty ]$.

\item[3'.] If $g_1$ is another Riemannian metric on $Q$ with sectional curvatures $\leq c <0$, and it is quasi-isometric to $g$
then the space at infinity
and the compactification of $Q$ are the same (as topological spaces) if taken with respect to $g_1$ or $g$.
(Note that we do not need $S$ to be totally geodesic with respect to $g_1$.)
\end{enumerate}

\vspace{.9in}

\noindent {\bf \Large  Section 2. Proof of Theorems A and B.}\\

We will say that two Riemannian metrics $g_0$, $g_1$ on a manifold $M$ are {\it homotopic } (or {\it isotopic})  if
there is a homotopy (or  isotopy)  $h_t:M\ra M$, $h_0=1_M$ such that $(h_1)^*g_0=g_1$.
We will also use the notation given at the beginning of section 1.
In what follows $M$ will denote a complete hyperbolic manifold with $dim\,\, M=n\geq 5$. 
The given hyperbolic metric will be denoted by
$g_0$. In what follows in this section if $M$ is non compact it is understood that all space of metrics considered are with respect to
$g_0$. For instance $\met(M)$ and $\mo(M)$ mean $\met(M,g_0)$ and $\mo(M,g_0)$ respectively. 
Furthermore $DIFF(M)=DIFF(M,g_0)$  and $DIFF_0(M)=DIFF_0(M,g_0)$.
Also, in this non compact case,
$g$, $g'$ in $\met(M)=\met(M,g_0)$ being homotopic means ``boundedly homotopic", that is, the homotopy $h_t$ is such that all
$h_t$ are at bounded $g_0$-distance from the identity (the bound independent of $t$).\\
\vspace{.6in}

\noindent {\bf \large 2.1. The map $\Lambda_\gamma^{\, \epsilon}:\met_\gamma^{\, \epsilon}(M)\ra P(\bS^1\times\bS^{n-2})$.}\\

Let's assume that there is an embedded closed 
geodesic $\gamma:\bS^1(\ell)\ra M$  in $M$ of length $\ell$, with orientable (hence trivial) normal bundle.
We define the subspace $\met _\gamma ^\epsilon(M)\sbs \mo(M)$ as the space of all metrics
$g\in\mo(M)$ such that:
\begin{enumerate}
\item[{\bf 1.}] The closed geodesic in $(M,g)$ representing the homotopy class of $\gamma$ coincides as a set with $\gamma$.
Moreover, the identity map $(\gamma ,g_0|_{\gamma})\ra (\gamma ,g|_{\gamma}))$ is a homothety
i.e. there is $c>0$ with $g_0(v,v)=cg(v,v)$, for all $v\in T\gamma$.

\item[{\bf 2.}]  $\perp^g_z\gamma= \perp^{g_0}_z\gamma$, for all $z\in \gamma$. 

\item[{\bf 3.}]  $g(u,v)=g_0(u,v)$, for all $u,v\in \perp^g_z\gamma$ and $z\in \gamma$.

\item[{\bf 4.}]  $exp^{\, g}_z(v)=exp^{\, g_0}_z(v)$, for all
$z\in \gamma$, $v\in \perp^g_z\gamma$ with $g(v,v)\leq\epsilon$.
\end{enumerate}

Define $ \met_\gamma(M)=\bigcup_{\epsilon >0}\met_\gamma^\epsilon (M)$ and define $\cT_\gamma(M)$ to be the image of $\met_\gamma(M)$ by
the quotient map $\met(M)\ra\cT(M)$.\\

Recall that $\R^+\times DIFF_0(M)$ acts on $\mo(M)$,
where $\R^+$ acts by scalar multiplication.
Let $\cD_\gamma(M)$ be the isotropy group of $\met_\gamma(M)$, that is: 

$$\cD_\gamma(M)=\Bigg\{ (\lambda, \phi)\in \R^+\times DIFF_0(M)
\, :\, \lambda\phi \Big(\met_\gamma(M)\Big)=\met_\gamma(M)\Bigg\}$$\\

\noindent {\bf Lemma 2.1.1.} {\it Let $(\lambda, \phi)\in \R^+\times DIFF_0(M)$. The following statements are equivalent.}
\begin{enumerate}
\item[{1.}] $(\lambda, \phi)\in \cD_\gamma(M)$.

\item[{2.}] {\it $\lambda\phi g\in\met_\gamma(M)$, for some $g\in\met_\gamma(M)$}.
 
\item[{3.}] {\it $\phi(\gamma)=\gamma$, the derivative $\sqrt{\lambda}\, D\phi_z:(T_zM,g_0)\ra (T_{\phi(z)}M,g_0)$ is an isometry, for all $z\in\gamma$,
and there is an $\epsilon >0$ such that $\phi(exp^{g_0}(v))=exp^{g_0}( D\phi(v))$, for all $v\in\perp^{g_0}\gamma$, $g(v,v) \leq \epsilon$.}
\end{enumerate}

\noindent {\bf Proof.} Clearly 1 implies 2. Also, an inspection of items 1-4 in the definition of $\met_\gamma(M)$ above
shows that 2 implies 3, and 3 implies 1. This proves the Lemma.\\

\noindent {\bf Lemma 2.1.2.} {\it Assuming $M$ is closed then the map $\met_\gamma(M)\ra\To_\gamma(M)$ 
is a  principal $\cD_\gamma(M)$-bundle.}\\

\noindent {\bf Proof.} Since the action of $DIFF_0(M)$ on $\mo(M)$ is free and $M$ is closed we have that the 
action of \,$\R^+\times DIFF_0(M)$\, on $\mo(M)$ is also free. This, together with
Ebin's Slice Theorem \cite{Ebin} implies that
$q:\mo(M)\ra\To(M)$ is a  principal $\Big(\R^+\times DIFF_0(M)\Big)$-bundle. Let $g\in \met_\gamma(M)$, then $g\in \met_\gamma^\delta(M)$,
for some $\delta >0$. It follows from Remark 2 after the proof of Prop. 1.1.1 that there is a $C^2$-open neighborhood
$W$ of $g$ in $\met(M)$ and a continuous map $r:W\ra \met_\gamma^\epsilon(M)$, for some $\epsilon >0$, such
that $r(g')$ is isotopic to $g'$, for every $g'\in W$, that is, $r(g')=\phi^* g'$, for some $\phi$ isotopic
to the identity. Hence $qr=q$.\\

$\Big[$Here to be able to apply  Remark 2 we have to use the fact,
due to Sampson \cite{Sampson} and Eells-Lemaire (\cite{EL2}, Prop. 5.5), that the the map $g\mapsto \eta_g$ is continuous,
where $\eta_g:\bS^1\ra M$ is the $g$-geodesic freely homotopic to $\gamma$.$\Big]$\\

Since $q:\mo(M)\ra\To(M)$ is a locally trivial $\Big(\R^+\times DIFF_0(M)\Big)$-bundle we can choose $W$ to be
a local product, i.e. there is an open set $V$ in $\To(M)$, an open neighborhood $U$ of $(1, 1_M)$ in $\Big(\R^+\times DIFF_0(M)\Big)$
and a map (a section)  $s:V\ra \mo(M)$ with $qs=1_V$ and
$W=U.s(V)=\{ \lambda\phi^* s(a)\, :\, a\in V\, ,\, (\lambda,\phi)\in U\}$. Note  that $q^{-1}(V)=\Big(\R^+\times DIFF_0(M)\Big). s(V)$.
Then $rs: V\ra\mo(M)$ is also a section (i.e. $qrs=qs=1_V$) and note that
the image of $rs$ lies in $\met_\gamma(M)$. It is straightforward to verify that 
$q|_{\met_\gamma(M)}^{-1}(V)=q^{-1}(V)\cap\met_\gamma(M)=\cD_\gamma(M).rs(V)$.
This together with Lemma 2.1.1 and the fact that $q$ is principal $\Big(\R^+\times DIFF_0\Big)$-bundle imply 
that the map $V\times\cD_\gamma(M)\ra q|_{\met_\gamma(M)}^{-1}(V)$ given by $(v,(\lambda, \varphi))\mapsto \lambda\varphi s(v)$
is a homeomorphism. This completes the proof of the Lemma.
\vspace{.4in}

We shall now define a map: $$\Lambda_\gamma^\epsilon:\met_\gamma^\epsilon(M)\ra P(\bS^1\times\bS^{n-2})$$
\noindent  Let $g\in\met^\epsilon_\gamma(M)$. Thus $g$ satisfies items 1-4 above.
Write also $N=N_\epsilon=\{ exp^{\, g_0}(v)\,\,:\,\, v\in\,\perp\gamma,\,\, g_0(v,v)\leq \epsilon\}$.
We have that $\perp^{\, g_0}\gamma=\perp^{\, g}\gamma$ and we just write $\perp \gamma$.
Let $Q$ be the covering space of $M$ corresponding to the infinite cyclic subgroup
of $\pi_1 (M, \gamma (\frac{\ell}{2\pi},0))$ generated by $\gamma$. 
(Here $(\frac{\ell}{2\pi} ,0)\in \bS^1(\ell)\sbs\R^2$.)
Denote also by $g_0$ and $g$ the pullbacks to $Q$
of the hyperbolic metric $g_0$ and the metric $g$.
Note that $\gamma$ and $N$ lift to $Q$ and we denote these liftings
also by $\gamma$ and $N$, respectively. Hence 1-4 above also hold true if we replace $M$ by $Q$.\\

We have that the normal exponential map
$exp^{\, g_0}:\perp\gamma\ra Q$ is a diffeomorphism, and since we are assuming that the normal
bundle of $\gamma$ is orientable (hence trivial) we have that 
$Q$ is diffeomorphic to $\bS^1\times\R^{n-1}$.
Therefore we will identify the following objects:
\begin{enumerate}
\item[$\bullet$] Identify $(Q,g_0)$ with $(\bS^1\times\R^{n-1}, \rho_0 )$
using one of the identifications given in the remark before 1.2.1, section 1.

\item[$\bullet$] Identify $\gamma\sbs Q=\bS^1\times\R^{n-1}$ with $\bS^1=\bS^1\times\{ 0\}\sbs \bS^1\times\R^{n-1}$.

\item[$\bullet$] Identify $\perp \gamma=\perp\bS^1$ also with $\bS^1\times\R^{n-1}$. Hence the exponential map
$exp^{\, g_0}$ is just the identity.

\item[$\bullet$] With all these identifications we have that $N=\bS^1\times\bS^{n-2}\times [0,\epsilon ]$.
\end{enumerate}

Then 1-4 above (with $Q$ instead of $M$) can be written in the following way:
\begin{enumerate}
\item[{\bf 1'.}] The unique closed geodesic in $(\bS^1\times\R^{n-1},g)$ (representing the homotopy class of $\bS^1$) coincides with $\bS^1$.
Moreover, $\bS^1(\ell_g )\ra (\bS^1, g|_{\bS^1})$ is an isometry, where $\ell_g=length_g(\bS^1)$.

\item[{\bf 2'.}]  $\perp^g_z\bS^1= \{ z\}\times\R^{n-1}$, for all $z\in \gamma$. 

\item[{\bf 3'.}]  $g(u,v)=\langle u,v\rangle_{\R^{n-1}}$, for all $u,v\in \R^{n-1}$ and $z\in \bS^1$.

\item[{\bf 4'.}]   $exp^{\, g}_z(v)=(z,v)$, for all
$(z,v)\in N$ with $\langle v,v\rangle_{\R^{n-1}}\leq\epsilon$. 
\end{enumerate}
Define $\varphi'_g\in DIFF(\bS^1\times\bS^{n-2}\times [\epsilon,\infty))$ as $\varphi'=exp^g$. Fix a diffeomorphism 
$\lambda:[\epsilon,\infty)\ra [0,1)$ and with it identify these two intervals to obtain $\varphi_g\in DIFF(\bS^1\times\bS^{n-2}\times [0,1))$.
(That is, $\varphi_g (z,v,t)= (z',v', t,)$ where  $(z', v',t')=exp^g(z,\lambda^{-1}(t)v)$ and $t=\lambda(t')$.)
By (4') we have $\varphi_g(z,v,0)=(z,v,0)$, hence $\varphi_g\in (DIFF(\bS^1\times\bS^{n-2}\times [0,1),\p )$.
We now extend $\varphi_g$ to $DIFF(\bS^1\times\bS^{n-2}\times [0,1])$ using asymptotics. We do this in the 
following way (for details see Section 2 of \cite{FO6}).
Take $(z,v)\in\bS^1\times\bS^{n-2}$. Then $\beta(t)=exp^g_z(t v)$, $t>0$, is a $g$-geodesic ray in $\bS^1\times\R^{n-1}$.
Hence it is a $\rho_0$-quasi-geodesic ray.
Therefore it determines a point at infinity $(\bar{z},\bar{v})\in\p_\infty(\bS^1\times\R^{n-1})=\bS^1\times\bS^{n-2}$.
(Equivalently, the Hausdorff $\rho_0$-distance between the $g$-geodesic $exp^g_z(tv)$ and the $\rho_0$-geodesic ray
$(\bar{z},t\bar{v})$, $t\geq 0$, is finite.)
We define then $\varphi_g(z,v,1)=(\bar{z},\bar{v},1)$.
An argument similar to the one given in the proof of Claim 3 of Section 3 in \cite{FO6} shows
that $\varphi_g:(\bS^1\times\bS^{n-2})\times I\ra (\bS^1\times\bS^{n-2})\times I$ is continuous
(this can also be proved using item 2' of Section 1.5).
So, we get $\varphi_g \in P(\bS^1\times\bS^{n-2})$.
We define then  $\Lambda_\gamma^\epsilon (g)=\varphi_g\in P(\bS^1\times\bS^{n-2})$.\\

\noindent {\bf Remark.} Taking $M=Q$ we also obtain a map $\Lambda^\epsilon_\gamma:\met^\epsilon_\gamma(Q)\ra P(\bS^1\times\bS^{n-2})$,
which is essentially the map $g\mapsto exp^g$, where $exp^g$ is the normal exponential map.
And we get the following commutative diagram (see Section 1.4):

$$\begin{array}{ccccc}\met_\gamma^\epsilon(M) & &\stackrel{\Lambda^\epsilon_\gamma}{ \longrightarrow}&&P(\bS^1\times\bS^{n-2})\\
\\
{\mbox{lift}}\downarrow &&\stackrel{{\scriptstyle{ \Lambda^\epsilon_\gamma}}}{}\nearrow \\
\\
\met^\epsilon_\gamma(Q)
\end{array}$$\\

\noindent {\bf Lemma 2.1.3} 
{\it If $M$ is closed the map $\Lambda^\epsilon_\gamma:\met^\epsilon_\gamma(M)\ra P(\bS^1\times\bS^{n-2})$ is continuous.}\\

\noindent{\bf Proof.} 
It is enough to prove that the map

$$\begin{array}{cccc}
F\,\,\,\,\,\,\,: &\met^\epsilon_\gamma(M)\times\Big( (\bS^1\times\bS^{n-2})\times [\epsilon, \infty ] \Big)& \longrightarrow & 
(\bS^1\times\bS^{n-2})\times [\epsilon, \infty ]\\
\\
& \Big( g, (z,v), t \Big)& \longmapsto & exp^g_z(tv)
\end{array}$$\\

\noindent is continuous. Choose $q=(\bar{z},\bar{v},\bar{t})\in (\bS^1\times\bS^{n-2})\times (\epsilon, \infty ]$
and a neighborhood of $q$ of the form $\bar{Z}\times\bar{V}\times (\bar{T},T]$ (we can have $T=\infty$).
Let also $F(g,z,v,t)=q$.
If $\bar{t}<\infty$ we can clearly find an open neighborhood $W$ of $(z,v)$ in $\bS^1\times\bS^{n-2}$
and a $C^2$-neighborhood $U$ of $g$ in $\met(M)$ such that $F(g',z',v')\in \bar{Z}\times\bar{V}\times (\bar{T},T]$,
for all $(z',v')\in W$ and $g'\in U\cap \met^\epsilon_\gamma(M)$. \\

Let $\bar{t}=\infty$. Then also $t=\infty$ and $T=\infty$. 
Since $t\mapsto exp^g_z(tv)$ is a $g$-geodesic ray, it is a  $(\lambda,\delta )$ $\rho_0$-quasi-geodesic ray, for some
$\lambda > 0$ and $ \delta \geq 0$.
Now, given $T>0$ we can find
an open neighborhood $W$ of $(z,v)$ in $\bS^1\times\bS^{n-2}$
and a $C^2$-neighborhood $U$ of $g$ in $\met(M)$ such that:

\begin{enumerate}
\item[{\it a.}] the $\rho_0$-distance between $exp^g_z(tv)$ and $exp^{g'}_{z'}(tv')$ is less than, say, one, for all $t<T$, 
$(z',v')\in W$ and $g'\in U$.

\item[{\it b.}] the $g'$-geodesic $t\mapsto exp^{g'}_{z'}(tv')$ is a $(2\lambda,\delta+1 )$ $\rho_0$-quasi-geodesic ray, for
every $g'\in U\cap\met^\epsilon_\gamma(M)$ and $(z',v')\in W$.

\end{enumerate}

But item 2' of Section 1.5 allows us to choose $T$ large enough so that
we can ensure that $exp^{g'}_{z'}(tv')\in \bar{Z}\times\bar{V}\times (\bar{T},\infty]$, for all
$(z',v', t)\in W\times [T,\infty]$ and  $g'\in U\cap\met^\epsilon_\gamma(M)$.
(Item 3' of section 1.5 is also used here because 2' refers to the topology generated by $g$, not $\rho_0$. But these two
topologies coincide, by item 3'.)
This proves the Lemma.
\vspace{.3in}

\noindent {\bf Addendum to Lemma 2.1.3} 
{\it If $M=Q$ then the map
$\Lambda^\epsilon_\gamma :\met^\epsilon_\gamma(Q)\ra P(\bS^1\times\bS^{n-2})$ is continuous.}\\\\

The proof of the Addendum is the same as the proof of lemma 2.1.3. Just recall that by $\met(Q)$ here we mean
$\met(Q,g_0)$ (or $\met(Q,\rho_0)$). Hence the $C^2$-neighborhood $U$ is really a $C^2$-neighborhood in (the true)
$\met(Q)$ intersected with a quasi-geodesic ball $B_{\lambda ,\delta}(g)$ (see section 1.4).\\

Now, note that the definition of $\Lambda_\gamma^\epsilon(g)=\varphi_g$ depends on $\epsilon$ because 
$\lambda:[\epsilon,\infty)\ra [0,1)$
depends on $\epsilon$. But $\tau\varphi_g$ does not depend on $\epsilon$, hence we get a well defined map
$\Lambda_\gamma:\met_\gamma (M)\ra TOP(\bS^1\times\bS^{n-2})$ defined by $\Lambda_\gamma (g)=\tau \Lambda_\gamma^\epsilon (g)
=\tau\varphi_g$.\\

\noindent {\bf Lemma 2.1.4} 
{\it If $M$ is closed, or $M=Q$, then $\Lambda_\gamma:\met_\gamma(M)\ra TOP(\bS^1\times\bS^{n-2})$ is continuous.}\\

\noindent{\bf Proof.} This follows from the fact that in the proof of  Lemma 2.1.3, the
$C^2$-neighborhood $U$ of $g$ in $\met(M)$ does not depend on $\epsilon$.
This proves Lemma 2.1.4.\\

Let $j>0$ be an integer. The $j$-sheeted cover $\bS^1\times\bS^{n-2}\ra \bS^1\times\bS^{n-2}$,
$(z,u)\mapsto (z^j,u)$ induces a continuous map $\nu_j:P^s(\bS^1\times\bS^{n-2})\ra P^s(\bS^1\times\bS^{n-2})$
obtained simply by pulling back (lifting) smooth pseudoisotopies using the $j$-sheeted cover.
Let $K\sbs P^s(\bS^1\times\bS^{n-2})$ be a compact subset and write $\iota$ for the inclusion of $K$ in 
$P^s(\bS^1\times\bS^{n-2})\sbs P(\bS^1\times\bS^{n-2})$. 
Also, let $\alpha:S\ra K$ be a map, with $S$ compact.
Let $r>0$ be given by Lemma 1.2.2. Hence $\Omega_\gamma(K)\sbs\mo(M)$.
We will also assume that $r>\epsilon$. Hence, by the definition of $\Omega_\gamma$ (see section 1.2) we have 
$$\Omega_\gamma(K)\sbs\met_\gamma^\epsilon(M)\sbs\met_\gamma(M)\sbs \mo(M)$$ Consider the following diagram:\\

$$\begin{array}{ccccccc}
S&\stackrel{\alpha}{\longrightarrow}& K& \stackrel{\Omega_\gamma}{\longrightarrow} & \met^\epsilon_\gamma(M)& 
\stackrel{i}{\hookrightarrow}&\met_\gamma(M)
\\ \\&&&\iota\searrow& \Lambda^\epsilon_\gamma\downarrow&&\Lambda_\gamma\downarrow\,\,\,\,\,\,\,\,\, \\ \\&& && P(\bS^1\times\bS^{n-2})&
\stackrel{\tau}{\longrightarrow}&TOP(\bS^1\times\bS^{n-2})
\end{array}$$\\

\noindent where $i$ denotes the inclusion. The square on the right is commutative, by the definition of $\Lambda_\gamma$. 
Let $c:S\ra P^s(\bS^1\times\bS^{n-2})$ denote the constant map 
$c(\varphi)=1_{\bS^1\times\bS^{n-2}\times [0,1]}$.\\

\noindent {\bf Proposition 2.1.5} 
{\it 
If, for some integer $j>0$, we have
$\nu_j\,\alpha\simeq c$  then $\iota\,\alpha\simeq\Lambda_\gamma^\epsilon\Omega_\gamma\,\alpha$, provided $r$ is large enough
(how large depending on $\ell$, $K$, $n$ and the homotopy between  $\nu_j\,\alpha $ and $ c$). Therefore, the following diagram
homotopy commutes}

$$\begin{array}{ccc}
S&\stackrel{\iota\,\alpha}{\longrightarrow}& P(\bS^1\times\bS^{n-2})\\ \\
{\scriptstyle{i\,\Omega_\gamma\,\alpha}}\downarrow& & \downarrow{\scriptstyle{\tau}}\\ \\
\met_\gamma(M)& \stackrel{\Lambda_\gamma}{\longrightarrow} &TOP(\bS^1\times\bS^{n-2})
\end{array}$$\\

\noindent The proof of this Proposition is given in Section 6.
\vspace{.8in}

\noindent {\bf \large 2.2. The map $\Delta_\gamma  :\To(M)^\bullet\ra  TOP(\bS^1\times\bS^{n-2})//L$.}\\

Recall from the introduction that $L\sbs TOP(\bS^1\times\bS^{n-2})$
is the subgroup of all ``orthogonal" self homeomorphisms of $\bS^1\times\bS^{n-2}$. That is, $f:\bS^1\times\bS^{n-2}\ra\bS^1\times\bS^{n-2}$
belongs to $L$ if $f(z,u)=(e^{i\theta}z, A(z)u)$, for some $e^{i\theta}\in \bS^1$, and $A:\bS^1\ra SO(n-1)$.
Also $\tau: P(\bS^1\times\bS^{n-2})\ra TOP(\bS^1\times\bS^{n-2})$ denoted the ``take top" map.\\

\noindent {\bf Proposition 2.2.1.} {\it Let $M$ be closed, or $M=Q$, and $g_1, g_2\in\met_\gamma(M)$ be homotopic. Then $ \Lambda_\gamma (g_1)=
\Lambda_\gamma (g_2)\, f$, for some $f\in L$.}\\

\noindent {\bf Remark.} Recall that for $M=Q$ we have that ``homotopic" means ``boundedly homotopic".\\

\noindent {\bf Proof.} Let $h:(M,g_1)\ra(M,g_2)$ be an isometry homotopic to the identity $1_M$.
Lifting $h$ to $\bS^1\times\R^{n-1}$ we obtain an isometry $h$ (we use the same letter) between
$(\bS^1\times\R^{n-1},g_1)$ and $(\bS^1\times\R^{n-1},g_2)$ such that $h$ is at bounded $\rho$-distance 
from the identity $1_{\bS^1\times\R^{n-1}}$.
(Since $g_1, g_2, \rho$ are all quasi-isometric, the same is true if we use the $g_i$-distance.) 
Let $Dh$ be the derivative of $h$. Since $g_1, g_2\in\met_\gamma(M)$ we have that:
\begin{enumerate}
\item[(i)] $h(\gamma)=\gamma$. Moreover, there is $e^{i\theta}\in\bS^1$ such that $h(z)=e^{i\theta}z$, for $z\in\bS^1$.

\item[(ii)] 
$h ( exp^{\, g_1} (  z, v)   )=
 exp^{\, g_2} ( e^{i\theta}z, Dh(z)v)$, for all 
$(z,v)\in\bS^1\times\R^{n-1}$.

\item[(iii)] $Dh(z):\{z\}\times\R^{n-1}\ra\{e^{i\theta}z\}\times\R^{n-1}$ is orthogonal with respect to $\rho$, $g_1$ and
$g_2$. (Recall that all these metrics coincide on $\perp\gamma=\bS^1\times\R^{n-1}$, see items 2 and 3 in the
definition of $\met_\gamma(M)$.) Therefore $Dh:\bS^1\times\bS^{n-2}\ra\bS^1\times\bS^{n-2}$ and $Dh\in L$.
\end{enumerate}

\noindent {\bf Claim.} {\it $\Lambda_\gamma(g_1)(z,v)=\Lambda_\gamma(g_2)(e^{i\theta}z,Dh(z)v)$.}\\

\noindent{\bf Proof of the Claim}. By item (ii) above and the fact that $h$ is at bounded distance from the
identity we have that the quasi-geodesics $t\mapsto exp^{\, g_1} (  z, tv)   $ and $t\mapsto exp^{\, g_2} ( e^{i\theta}z, tDh(z)v)$
are a bounded distance apart. Hence, by the definition of $\Lambda_\gamma(g)(z,v,1)$ we have
$\Lambda_\gamma(g_1)(z,v,1)=\Lambda_\gamma(g_2)(e^{i\theta}z,Dh(z)v,1)$. The claim follows.\\

If $f:\bS^1\times\bS^{n-2}\ra\bS^1\times\bS^{n-2}$ is defined as $f(z,v)=(e^{i\theta}z, Dh(z)v)$, then
from the claim we have that
$\Lambda_\gamma (g_1)= \Lambda_\gamma (g_2)f$. This proves the Proposition.\\

It follows from the proof of Proposition 2.2.1 that the function $f$ is just the derivative  $Dh$. By modifying the proof
of the Proposition in a straightforward way we obtain the following addition this Proposition:\\

\noindent {\bf Addendum to Proposition 2.2.1.} {\it Let $M$ be closed, or $M=Q$, and $g_1, g_2\in\met_\gamma(M)$ such that $\lambda\varphi g_1=g_2$,
for some $(\lambda,\varphi)\in \cD_\gamma(M)$. Then $ \Lambda_\gamma (g_1)=
\Lambda_\gamma (g_2)\, f$, where $f=\sqrt{\lambda}D(\varphi)\in L$.}\\

Therefore
we obtain a continuous group homomorphism
$D:\cD_\gamma(M)\ra TOP(\bS^1\times\bS^{n-2})$, $D(\lambda,\varphi)=\sqrt{\lambda}D(\varphi)$, 
such that $\Lambda_\gamma(\varphi g)=\Lambda_\gamma (g)D(\varphi)$, for all
$g\in\met_\gamma(M)$. This together with Lemma 2.1.2 and items 2 and 3 of Section 1.3 imply the following Proposition.\\

\noindent {\bf Proposition 2.2.2.} {\it The following diagram commutes, where the lower horizontal arrow 
is the `orbit map' induced by $\Lambda_\gamma$.}\\

$$\begin{array}{ccc}
  \met_\gamma (M)^\bullet& 
\stackrel{\Lambda_\gamma^\bullet}{\longrightarrow}&TOP(\bS^1\times\bS^{n-2})^\bullet\\
 \\ \downarrow&&\downarrow\,\,\,\,\,\,\,\,\, \\ \\ \met_\gamma(M)//\cD_\gamma(M)&
\longrightarrow&TOP(\bS^1\times\bS^{n-2})//L
\end{array}$$\\\\

\noindent {\bf Remark.} If $M$ is closed item 3 of Section 1.3 together with Lemma 2.1.2 give a canonical identification of
$\To_\gamma(M)^\bullet$ with $\met_\gamma(M)//\cD_\gamma(M)$.\\

\noindent {\bf Lemma 2.2.3.} {\it We have that }$$\met_\gamma(Q)//\cD_\gamma(Q)=\mo (Q)//\R^+\times DIFF_0(Q)$$

\noindent {\bf Proof.} Since $\cD_\gamma(Q)$ is the isotropy group of $\met_\gamma(Q)$ we have that the semi-simplicial map
$\cS(\met_\gamma(Q))/\cS(\cD_\gamma(Q))\ra \cS(\mo (Q))/\cS(\R^+\times DIFF_0(Q))$ is injective.
Surjectivity follows from the following three facts:

\noindent {\bf 1.} Any closed $g$-geodesic representing $\gamma$ has to be embedded. Hence Proposition 1.1.1
implies that $\cT_\gamma(Q)//\cD_\gamma(Q)=\To(Q)$. 

\noindent {\bf 2.} $k$-simplices are contractible.

\noindent {\bf 3.} 
The fact,
due to Sampson \cite{Sampson} and Eells-Lemaire (\cite{EL2}, Prop. 5.5), that the the map $g\mapsto \eta_g$ is continuous,
where $\eta_g:\bS^1\ra M$ is the $g$-geodesic freely homotopic to $\gamma$.
This proves the Lemma.\\

Now, let $M$ be closed. Lifting metrics gives us a map

{\scriptsize $$\mo (M)//\R^+\times DIFF_0(M)\ra\mo (Q)//\R^+\times DIFF_0(Q)$$}

\noindent But $\mo(M)\ra\To(M)$ is a (locally trivial) bundle, hence 
the domain of the above map is canonically identified with
$\To(M)^\bullet$ (see section 1.3). Denote by $\Delta_\gamma$ the composition map

{\scriptsize $$\To(M)^\bullet\ra \mo (Q)//\R^+\times DIFF_0(Q)=\met_\gamma(Q)//\cD_\gamma(Q)\ra TOP(\bS^1\times\bS^{n-2})//L$$}

\noindent and by chasing diagrams around we obtain the following commutative diagram

$$\begin{array}{ccc}
  \met_\gamma (M)^\bullet& 
\stackrel{\Lambda_\gamma^\bullet}{\longrightarrow}&TOP(\bS^1\times\bS^{n-2})^\bullet\\
 \\ \downarrow&&\downarrow\,\,\,\,\,\,\,\,\, \\ \\ \To(M)^\bullet&
\stackrel{\Delta_\gamma}{\longrightarrow}&TOP(\bS^1\times\bS^{n-2})//L
\end{array}$$\\

\noindent {\bf \large 2.3. Completion of the Proof of Theorems A and B.}\\

Let $n$ and $k$ be such that
the map $\Upsilon=\Upsilon _{n,\, k}:\pi_k (\, P^s(\bS^1\times\bS^{n-2})^\bullet\, )\ra \pi_k(\, TOP(\bS^1\times\bS^{n-2}) \, // \, L \,)$
is strongly nonzero. Let $\beta:\bS^k\ra P^s(\bS^1\times\bS^{n-2})^\bullet$ represent a class such that 
$\Upsilon([\beta ])\neq 0$. Using the fact that the map $P^s(\bS^1\times\bS^{n-2})^\bullet\ra P^s(\bS^1\times\bS^{n-2})$ is a
weak homotopy equivalence we can assume that $\beta$ has the form 
$\bS^k\longrightarrow (\bS^k)^\bullet\stackrel{\alpha^\bullet}{\longrightarrow}P^s(\bS^1\times\bS^{n-2})^\bullet$, 
for some $\alpha:\bS^k\ra P^s(\bS^1\times\bS^{n-2})$, and some (fixed) homotopy equivalence $\bS^k\ra (\bS^k)^\bullet$.
We intend to use Proposition 2.1.5, so we need to know that we can choose $\alpha$ so that 
there is an integer $j>0$ such that $\nu_j\,\alpha$ is nullhomotopic, provided  $r$ is large enough.
This follows from  Lemma 2.3.1 below, that shows, together with the definition of `strongly nonzero'
(see sentence before Theorem A), that in fact we have infinitely many choices
for $\alpha$.\\

\noindent {\bf Lemma 2.3.1.}
{\it Let $T$ be an infinite torsion subgroup of $\,\pi_k \Big(P^s (\bS^1\times\bS^{n-2})\Big)$
and assume that $n >> k$. Then the subgroup $\overline{T}$ of $T$ consisting of all elements
 which vanish under some $(\nu_j)_*$ is also infinite. } \\
 
 The proof of Lemma 2.3.1 is given in section 5.\\

 \noindent {\bf Remarks.} \\

\noindent 1. Recall that Theorem C, proved in section 3, shows that such subgroups $T$ do exist in our relevant cases.\\
2. That $\overline{T}$ is a subgroup follows easily from the fact that $\nu_j \nu_k =\nu_k \nu_j$ 
for all positive integers $j$ and $k$.\\
3. $n>>k$ in Lemma 2.3.1 and the rest of Section 2 refers to Igusa's stable range,
namely $n>$ max $\{ 3k+8,\, 2k+9\}$ 
(see \cite{I}, p.6).\\

Write $K=\alpha(\bS^k)$. Then $K$ is compact.
It follows from Lemma 2.3.1  above, Lemma 1.2.2  and Proposition 2.1.5  that we can choose $r>0$ large enough so that  
 $\Omega_\gamma(K)\sbs\mo(M)$ is well defined  and  the following diagram commutes, up to homotopy

$$\begin{array}{cccc}
\bS^k&\stackrel{\iota\,\alpha}{\longrightarrow}& P(\bS^1\times\bS^{n-2})\\ \\
{\scriptstyle{i\,\Omega_\gamma\,\alpha}}\downarrow& & \downarrow{\scriptstyle{\tau}}\\ \\
\met_\gamma(M)& \stackrel{\Lambda_\gamma}{\longrightarrow} &TOP(\bS^1\times\bS^{n-2})
\end{array}$$\\

\noindent Therefore, the upper right square of the following diagram homotopy commutes:\\

$$\begin{array}{ccccc}

&&\bS^k& \stackrel{\iota^\bullet\,\beta}{\longrightarrow} &   P(\bS^1\times\bS^{n-2})^\bullet   
\\ \\
&&{\scriptstyle{(i\,\Omega_\gamma)^\bullet\,\beta}}\downarrow\,\,\,\,\,\,\,\,\,\,\,\,\,\,\,\,\,\,& &{\scriptstyle{\tau^\bullet}}\downarrow\,\,\,\,\,\,\,
\\ \\ 
\met_\gamma (M)&\leftarrow& \met_\gamma(M)^\bullet&\stackrel{\Lambda_\gamma^\bullet}{\longrightarrow}&TOP(\bS^1\times\bS^{n-2})^\bullet
\\ \\ 
\downarrow&& \downarrow&&\downarrow 
\\ \\
\To(M)&\leftarrow& \To(M)^\bullet&\stackrel{\Delta_\gamma}{\longrightarrow}&TOP(\bS^1\times\bS^{n-2})//L

\end{array}$$\\

\noindent The left square obviously commutes and the bottom right square is the commutative diagram at the end of section 2.2.\\

 Hence, applying the $\pi_k$ functor to the diagram above we obtain a commutative diagram of groups,
 and, using the facts that $\Upsilon([\beta])\neq 0$, and $\pi_k(X^\bullet)\ra \pi_k(X) $ is an isomorphism for
 any $X$, we get that the map $\pi_k(\met_\gamma(M))\ra\pi_k(\To(M))$ is nonzero. But this map
 factors through $\pi_k(\mo(M)))\ra\pi_k(\To(M))$. This proves Theorem A.\\
 
 To prove Theorem B  apply the functor $H_k$ instead of $\pi_k$ to the above diagram and use the fact that $h\Upsilon([\beta])\neq 0$. 
The rest of the proof is similar. This proves Theorem B.
\vspace{.8in}

\noindent {\bf \Large  Section 3. Proof of Theorem C.}\\

Throughout this and next Sections we will use the following notation. For an abelian group $A$,
$\tau(A)$ is the torsion subgroup of $A$. 
Also, for a prime $p$,  $\tau_p(A)$ is the subgroup of $A$ consisting of all elements of order a power
of $p$. Then $\tau(A)=\bigoplus_{p\,\, prime}\tau_p(A)$ and $A/\tau(A)$ is torsion free.\\

For a prime $p$, let $\zp$ be the ring $\Z$ localized at $p$, i.e.
$\zp=\Z[\frac{1}{2}, \frac{1}{3},...,\hat{\frac{1}{p}},...]=\{\frac{r}{s}\in\Q\,\,:\,\, (r,s)=1,\,\, p \not\,\mid \, s\, \}$.
We denote by $\cC_p$ the class of all abelian groups $A$ for which 
$A\otimes \zp$ is finitely generated as a $\zp$-module. We remark that if $A$ is in $\cC_p$ then 
$\tau_p(A)$ is finitely generated, hence finite (see  Section 4).
Also, recall that if $\cA$ is a class of abelian groups, a group
homomorphism $f:G_1\ra G_2$ is an $\cA$-isomorphism if $ker f, coker f$ are in $\cA$.\\

Theorem C for the case $k=2p-4$ is a direct consequence of Theorems D and E below. 
Theorem D is a version of Theorem C for $\Upsilon$.
In both Theorems D and E (and in the paragraph following the statement of Theorem E)
the base point for the homotopy groups $\pi_i(\,\,\, )$, $i>0$, is the one corresponding to the identity
homeomorphism.\\

\noindent {\bf Theorem D.} {\it Consider the map $\Upsilon_{n,\, k}$. We have the following cases:}
\begin{enumerate}
\item[{}] {\bf k=0}\newline
{\it the group $\pi_0 ( P^s(\bS^1\times\bS^{n-2}))$ contains a subgroup $(\Z_2)^\infty$
and $\Upsilon_{n,\, 0}$ restricted to this  $(\Z_2)^\infty$ is injective, provided $n\geq 10$.}

\item[{}] {\bf k=1}\newline
{\it the group $\pi_1 ( P^s(\bS^1\times\bS^{n-2}))$ contains a subgroup $(\Z_2)^\infty$
and $\Upsilon_{n,\, 1}$ restricted to this  $(\Z_2)^\infty$ is injective, provided $n\geq 12$.}

\item[{}] 
{\bf k=2p-4, p$>$2 prime.}\newline
{\it the group $\pi_k ( P^s(\bS^1\times\bS^{n-2}))$ contains a subgroup $(\Z_p)^\infty$
and $\Upsilon_{n,\, k}$ restricted to this  $(\Z_p)^\infty$ is injective, provided $n\geq 3k+8$.}

\end{enumerate}

\vspace{.3in}

\noindent {\bf Theorem E.} {\it For $p>2$ prime, and $n>3k+8$, the Hurewicz map
$h:\pi_{2p-4}(TOP(\bS^1\times\bS^{n-2})//L)\ra H_{2p-4}(TOP(\bS^1\times\bS^{n-2})//L)$ 
is a $\cC_p$-isomorphism.}\\

Note that Theorem C 
for the case $k=0$, follows directly from
Theorem D above (case $k=0$) because $h:\pi_0(TOP(\bS^1\times\bS^{n-2})//L)\ra H_0(TOP(\bS^1\times\bS^{n-2})//L)$ 
is one-to-one. The case $k=1$ of Theorem C follows also from Theorem D (case $k=1$) because 
$\pi_1(TOP(\bS^1\times\bS^{n-2})//L)$ is abelian (therefore 
$h:\pi_1(TOP(\bS^1\times\bS^{n-2})//L)\ra H_1(TOP(\bS^1\times\bS^{n-2})//L)$ is an isomorphism).
And $\pi_1(TOP(\bS^1\times\bS^{n-2})//L)$ is abelian because 
$\pi_1(TOP(\bS^1\times\bS^{n-2}))$ is abelian and
$\pi_1(TOP(\bS^1\times\bS^{n-2}))\ra \pi_1(TOP(\bS^1\times\bS^{n-2})//L)$  is onto
(for this last fact see the proof of Proposition 3.2 below, after the proof of Prop. 3.4).\\

\noindent {\bf Proof of Theorem D.} Write $N=\bS^1\times\bS^{n-2}$ and let $F:P^s(N)\ra P(N)$ be the ``forget structure map".
We first treat the case {\bf k=2p-4, p$>$2 prime.}
Note that $L$ is homeomorphic to $\bS^1\times SO(n-1)\times\Omega\Big( SO(n-1) \Big)$,
thus its homotopy groups are all finitely generated. In particular, the natural map
$\pi_k(TOP(N))=\pi_k(TOP(N)^\bullet)\ra\pi_k(TOP(N)//L)$ is a $\cC_p$-isomorphism.
Since $\Upsilon_{n,\,\, k}$ is the composition of $\pi_k(\tau)\circ\pi_k(F)$ with
the natural map,
to prove Theorem D for the case {\bf k=2p-4, p$>$2 prime}, it is enough to 
prove that there is a subgroup $(\Z_p)^\infty$ of $\pi_k(P^s(N))$ such that
$\pi_k(\tau)\circ\pi_k(F)$ restricted to $(\Z_p)^\infty$ is one-to-one.
We shall prove this.\\

There is a spectral sequence (see \cite{Hat}) with $E^1_{st}=\pi_t(P(N\times I^s))$
converging to $\pi_{s+t+1}\Big( \widetilde{TOP}(N) //TOP(N) \Big)$.\\

\noindent {\bf Remark.} 
With our notation we have  $\widetilde{TOP}(N) //TOP(N)=\Big|\widetilde{TOP}(N)/S\Big( TOP(N)\Big) \Big|$, 
where $\widetilde{TOP}(N)$ is the space of block topological automorphisms of $N$.
This is the simplicial set
whose $k$-simplices are the automorphisms of $N\times\Delta^k$ which leave invariant each
$N\times \Big( {\mbox{face of }}\Delta^k\Big)$. The space of block self homotopy equivalences
$\widetilde{G}(N)$ is defined analogously and we have that
$\widetilde{G}(N)\simeq G(N)$, where $G(N)$ is the H-space of self homotopy equivalences of $N$.  
The corresponding quotients are 
similarly defined (see \cite{Hat}).\\

In particular we have $E^1_{0\, t}=\pi_t(P(N))$. Consider the composite map\\

$$
\begin{array}{ccccccccc}
\pi_t(P(N))& = & E^1_{0\, t} & \stackrel{{\mbox{\tiny onto}} }{\longrightarrow} & E^\infty_{0\, t} 
& \stackrel{{\mbox{\tiny 1-1}}}{\longrightarrow} & \pi_{t+1}\Big( \widetilde{TOP}(N)//TOP(N) \Big) & \longrightarrow & \pi_t(TOP(N)) \\ 
 \\      
 & & {\mbox{\tiny onto}} \searrow & &  \nearrow {\mbox{\tiny onto}}       \\
 \\
 & &  & E^2_{0\, t} &     
\end{array}
$$
\\

\noindent This map can be identified with $\pi_t(\tau): \pi_t(P(N))\ra\pi_t(TOP(N))$.
In Igusa's stable range we can identify $E^2_{s\, t}$ as $H_s\Big(\Z_2;\pi_t \cP(N) \Big)$,
where $\cP(\,\,\, )$ is the stable pseudo-isotopy functor. Then $E^2_{s\,t}$ is a 2-torsion
group when $s> 0$ and $s+n\geq 3t+8$. And the Igusa stable condition holds for both
$E^2_{k\, 0}$ and $E^2_{s\, t}$ such that both $s+t=k+1$ and $s\geq 2$. Consequently
the surjective map $E^2_{0\, k}\ra E^\infty_{0\, k}$ is $\cC_p$-injective (i.e. its kernel is in
$\cC_p$). \\

\noindent {\bf Claim.} {\it The map $\pi_{t+1}\Big( \widetilde{TOP}(N)//TOP(N) \Big) \longrightarrow \pi_t(TOP(N))$
is a $\cC_p$-isomorphism.}\\

\noindent {\bf Proof of the Claim.} To see this first observe that $\pi_i(G(N))$
is finitely generated, hence we need only to show that $\pi_i\Big( \widetilde{G}(N)//\widetilde{TOP}(N)  \Big)$
is finitely generated. For this we use the functional space approach to surgery theory developed by
Quinn in his thesis \cite{Quinn} and exposed in \cite{W}, pp. 240-241. In particular, there is a fibration (up to homotopy)
$\widetilde{G}(N)//\widetilde{TOP}(N)\ra (G/TOP)^N\ra {\cal{L}}(N)$, where the homotopy groups of
${\cal{L}}(N)$ are the Wall surgery groups of $\pi_1(N)=\Z$ and hence finitely generated due to Browder \cite{Browder}.
The Homotopy groups of $(G/TOP)^N$ are finitely generated due to Kirby-Siebenmann
\cite{KiSi}. This proves the claim.\\

It follows from the Claim and the discussion above that the composite

$$ E^2_{0\, k}  \stackrel{{\mbox{\tiny onto}} }{\longrightarrow}  E^\infty_{0\, k} 
 \stackrel{{\mbox{\tiny 1-1}}}{\longrightarrow}  \pi_{k+1}\Big( \widetilde{TOP}(N)//TOP(N) \Big)  \longrightarrow  \pi_k(TOP(N)) $$

\noindent is $\cC_p$-injective. Now, 
the surjective map $\pi_k(\cP(N))=E^1_{0\, k}\ra E^2_{0\, k}$
can be identified as the quotient map $$\pi_k(\cP(N))\ra H_0\Big( \Z_2; \pi_k\cP(N) \Big)$$
\noindent where $\Z_2$ acts on $\cP(N)$ via the ``turning upside down" involution $-$ on $\cP(N)$.
Therefore, to prove Theorem D, 
for the case {\bf k=2p-4, p$>$2 prime}, it is enough to 
prove that there is a subgroup $(\Z_p)^\infty$ of $\pi_k(\cP^s(N))$ such that
the map
$$\pi_k(\cP^s(N))\stackrel{\pi_k(F)}{\longrightarrow}\pi_k(\cP(N))\ra H_0\Big( \Z_2; \pi_k\cP(N) \Big)$$
\noindent restricted to $(\Z_p)^\infty$ is one-to-one.
But observe that the map $\pi_k(F):\pi_k\cP^s(N)\ra\pi_k\cP(N)$ is a $\Z_2$-module map which is a $\cC_p$-isomorphism
(see Lemma 4.1 of \cite{FO6}). Consequently, the right hand vertical arrow in the following diagram is also a
$\cC_p$-isomorphism:

$$\begin{array}{ccccc}
\pi_k(\cP^s(N)) & & \longrightarrow & & H_0\Big( \Z_2\, ;\, \pi_k(\cP^s(N))  \Big)\\
\\
{\scriptstyle \pi_k(F)}\downarrow\,\,\,\,\,\,\,\,\,\,\,\,\,&& & &\downarrow \\
\\
\pi_k(\cP(N)) & & \longrightarrow & & H_0\Big( \Z_2\, ;\, \pi_k(\cP(N))  \Big)

\end{array}
$$

Consequently, to prove Theorem D, 
for the case {\bf k=2p-4, p$>$2 prime}, it is enough to 
prove that there is a subgroup $(\Z_p)^\infty$ of $\pi_k(\cP^s(N))$ such that
the map
$$\pi_k(\cP^s(N))\ra H_0\Big( \Z_2; \pi_k\cP^s(N) \Big)$$
\noindent restricted to $(\Z_p)^\infty$ is one-to-one.
We will prove this.\\

There is a $\Z_2$-module map $\pi_{k+2}(A(N))\ra\pi_k(\cP(N))$ which is both an epimorphism and
a $\cC_p$-isomorphism, where $A(\,\,\, )$ is Waldhausen's functor. (See section 4 of \cite{FO6}
for more details.) Therefore we have a commutative diagram

$$\begin{array}{ccccc}
\pi_{k+2}(A(N)) & & \longrightarrow & & H_0\Big( \Z_2\, ;\, \pi_{k+2}(A(N))  \Big)\\
\\
\downarrow&& & &\downarrow \\
\\
\pi_k(\cP^s(N)) & & \longrightarrow & & H_0\Big( \Z_2\, ;\, \pi_k(\cP^s(N))  \Big)

\end{array}
$$

\noindent such that the right hand vertical arrow is also a $\cC_p$-isomorphism.
But an obvious modification of the argument proving Prop. 4.6 of \cite{FO6} yields
a subgroup $(\Z_p)^\infty$ of $\pi_{k+2}(A(N))$ which maps monomorphically into
$H_0\Big( \Z_2\, ;\, \pi_{k+2}(A(N))  \Big)$, and therefore the same is true
for $\pi_k(\cP^s(N))$ and the map $\pi_k(\cP^s(N))\ra H_0\Big( \Z_2\, ;\, \pi_k(\cP^s(N))  \Big)$,
which is what we wanted to prove. This proves Theorem D for the case
{\bf k=2p-4, p$>$2 prime}.
\vspace{.3in}

It can readily be checked that for the case  {\bf k=0} no changes are needed
and the whole argument goes through, except for the part in which we prove that
the map $E^2_{0k}\stackrel{{\mbox {\tiny onto}}}{\longrightarrow} E^\infty_{0k}$ is $\cC_p$-injective. Since we are working with
prime {\bf p=2} in this case, we can not use the fact that the terms $E^2_{st}$ are 2-torsion,
for $s>0$. But note that $\pi_j(\cP(N))=0$ for $j<0$, hence the spectral sequence is a first
quadrant spectral sequence and it follows that $E^2_{00}=E^\infty_{00}$.\\

Similarly, for the case  {\bf k=1} the only problem appears in the proof of the  $\cC_p$-injectivity of
the map $E^2_{0k}\stackrel{{\mbox {\tiny onto}}}{\longrightarrow} E^\infty_{0k}$. 
Again, since this is a first quadrant spectral sequence, we have that $E^3_{01}=E^\infty_{01}$ and we 
obtain the following exact sequence
$$H_2\Big( \Z_2;\pi_0\cP(N) \Big)=E^2_{20}\ra E^2_{01}\ra E^3_{01}=E^\infty_{01}\ra 0\hspace{.4in}$$

But observe that $E^2_{20}=H_2\Big( \Z_2;\pi_0\cP(N) \Big)$ is $\cC_2$-isomorphic to
$H_2\Big( \Z_2;\pi_2 A(N) \Big)$, because the composite $\Z_2$-module map
$\pi_2A(N)\ra \pi_0 \cP^s(N)\ra\pi_0\cP(N)$ is a $\cC_2$-isomorphism. Furthermore, we have
(see discussion in the last Section of \cite{FO6})
$$\pi_2A(N)\cong \pi_2 A(\bS^1)=\pi_2(A(*))\oplus\pi_{1}(A(*))\oplus\pi_2(N_-A(*))\oplus\pi_2(N_+A(*))$$
\noindent and the conjugation leaves invariant the first two terms and interchanges the last two.
But $\pi_2(A(*))$ and $\pi_{1}(A(*))$ are both finitely generated, hence $H_2\Big( \Z_2;\pi_2 A(\bS^1) \Big)$
is finitely generated. Consequently $E^2_{20}=H_2\Big( \Z_2;\pi_0\cP(N) \Big)$ is in $\cC_2$. Therefore
$ E^2_{01}\ra E^3_{01}$ is $\cC_2$-injective. This concludes the proof of Theorem D.\\\\

To prove Theorem E we want to use the following general version of Hurewicz's Theorem (see Spanier \cite{Sp}, p. 510):\\

\noindent {\bf Theorem.} {\it Let $X$ be a strongly simple space and $\cA$ an acyclic Serre ring of abelian groups.
If $\pi_j(X)\in\cA$, $1\leq j<k$ then $H_j(X)\in \cA$, $1\leq j<k$ and the Hurewicz map
$h:\pi_k(X)\ra H_k(X)$ is an $\cA$-isomorphism.}\\

For the definition of an {\it acyclic Serre ring of abelian groups}\, and
{\it strongly simple space} see Spanier \cite{Sp}, Chap.9, Sec. 6. 
Using the general version of Hurewicz's Theorem given above, Theorem E reduces to the following three Propositions.\\

\noindent {\bf Proposition 3.1.} {\it The class $\,\cC_p$ is an acyclic Serre rings, for any prime $p$.}\\

\noindent {\bf Proposition 3.2.} {\it The space $TOP(\bS^1\times\bS^{n-2})//L$ is strongly simple.}\\

\noindent {\bf Proposition 3.3.} {\it For $p>2$ prime, we have that $\pi_j(TOP(\bS^1\times\bS^{n-2})//L)\in\cC_p$,
$1\leq j<2p-4$.}\\

The proof of Proposition 3.1 is given in Section 4. The proof of Proposition 3.3 is given at the end of this Section.
Before we prove Proposition 3.2,  we give first a somewhat more general result:\\

\noindent {\bf Proposition 3.4.} {\it Let $G$ be a topological group and  $H$ a subgroup of $G$.
Let $p:G^\bullet\ra G//H$ be the projection and assume that $\pi_1(p):\pi_1(G^\bullet )\ra\pi_1(G//L)$ is onto.
Then $G//H$ is strongly simple.} \\

\noindent {\bf Proof.} The identity of $G$ will be denoted by $e$ and the corresponding vertex in $G^\bullet$
will also be denoted by $e$.
According to Example 18 of \cite{Sp} (p. 510), it is enough to to prove the following:
for each $a\in \pi_1(G//H, p(e))$ there is a map $\omega_a:\bS^1\times G//H\ra G//H$
such that $\omega_a|_{\,\bS^1\times \{ p(e)\}}$ represents $a$ and $\omega_a|_{\{ 1\}\times G//H}$
is homotopic to the identity.\\

Let $a\in \pi_1(G//H, p(e))$. Since $\pi_1(p)$ is onto there is a loop $\alpha':\bS^1\ra G^\bullet$
such that $p\,\alpha'$ represents $a$. After composing $\alpha'$ with the projection map $q:G^\bullet\ra G$, we get
a loop $\alpha=q\alpha' :\bS^1\ra G$.\\

Identify $\bS^1$ with the boundary $\p\Delta^2$ of the canonical 2-simplex, with its canonical simplicial complex structure,
that is, the one with three vertices: $e_0$, $e_1$, $e_2$, and three 1-simplices $[e_0,e_1]$, $[e_1,e_2]$, $[e_0,e_2]$.
Let $\Sigma$ be the simplicial set induced by this structure. 
Note that
all $n$-simplices of $\Sigma$, $n>1$, are degenerate. Note also that the geometric realization $|\Sigma|$ is canonically homeomorphic to $\bS^1$
and we just write $|\Sigma|=\bS^1$.
The set of $n$-simplices that form $\Sigma_n$ are 
sequences of $n+1$ vertices of the form $e_i ...e_ie_j...e_j$, $i\leq j$, $i,j\in\{ 0,1,2\}$. 
Such an object is determined by three integers $i,j,k$ where $i,j$ are as before and $k$ is the number of times
$e_i$ appears in the sequence (hence $e_j$ appears $(n+1)-k$ times). We denote this $n$-simplex by $\tau^n(i,j,k)$. 
For each $n$-simplex $\tau=\tau^n(i,j,k)$ denote by $\bar{\tau}:\Delta^n\ra [e_i,e_j]$ the simplicial map
that sends the first $k$ vertices of $\Delta^n=[e_0,...,e_{n}]$ to $e_i$ and the last $(n+1)-k$ to $e_j$.\\

Consider the simplicial set $\Sigma\times S(G)$ and recall that $(\Sigma\times S(G))_n=\Sigma_n\times S(G)_n$. 
Since $\, |\Sigma|=\bS^1$ is a CW-complex we have that $|\Sigma\times S(G)|=|\Sigma|\times |S(G)|=\bS^1\times G^\bullet$
(see \cite{LW}, p. 97).
We now define a simplicial map $\Omega'=\Omega'_a:\Sigma\times S(G)\ra S(G)$ in the following way. For $(\tau, \sigma)\in
\Sigma_n\times S(G)_n=(\Sigma\times S(G))_n$,  define $\Omega' (\tau,\sigma):\Delta^n\ra G$
as $\Omega'(\tau,\sigma)(v)=\alpha (\bar{\tau}(v))\, .\,\sigma(v)$, $v\in\Delta^n$.
Applying the geometric realization functor we obtain a map $|\Omega'|=|\Omega'_a|:\bS^1\times G^\bullet\ra G^\bullet$.
Write $\omega'=\omega_a'=|\Omega'|$.\\

\noindent {\bf Claim 1.} {\it The map $\omega'$ restricted to $\bS^1\times\{ e\}$ represents $\alpha'$.}\\

\noindent{\bf Proof of Claim 1.} The inclusion $\iota: \{ e\}\hookrightarrow G$ induces the simplicial map $S(\iota):S(\{ e\} )\ra S(G)$.
Note that $S(\{ e\})$ has exactly one $n$-simplex: $\sigma_e^n\Delta^n\ra\{ e\}$. Consider the following sequence of simplicial maps:
$$\Sigma\ra\Sigma\times S(\{ e\} )\stackrel{S(1_\Sigma)\times S(\iota)}{\longrightarrow}\Sigma\times S(G)\stackrel{\Omega'}{\ra}S(G) $$

\noindent where the first map is given by $\tau\mapsto (\tau,\sigma_e^n)$.
It can be easily verified that the image of a $\tau\in\Sigma$ in $S(G)$ after applying this sequence
of simplicial maps is $\alpha\,\bar\tau$. Hence the image of $[e_i,e_j]\in\Sigma$ by this sequence of maps is
the singular 1-simplex $\Delta^1=[e_0,e_1]\ra[e_i,e_j]\stackrel{\alpha}{\ra} G$. It follows
that after applying the geometric realization functor to the sequence above and composing with $q:G^\bullet\ra G$
at the end we obtain 
$$\bS^1\ra\bS^1\times\{ e\} \stackrel{1_{\bS^1}\times i}{\longrightarrow}\bS^1\times G^\bullet\stackrel{\omega'}{\ra}G^\bullet 
\stackrel{q}{\ra}G$$
\noindent (here $i=|S(\iota)|$ is the inclusion) and this composition is just $\alpha$. Since $\pi_1(q):\pi_1(G^\bullet)\ra\pi_1(G)$
is an isomorphism, the claim follows.\\

\noindent {\bf Claim 2.} {\it The map $\omega'$ restricted to $\{ 1\}\times G^\bullet$ is the identity.}\\

\noindent{\bf Proof of Claim 2.} The proof is similar to the proof of claim 1. Just consider the sequence of obvious simplicial
maps
$$S(G)\ra\{ 1\} \times S(G )\ra\Sigma\times S(G)\stackrel{\Omega'}{\ra}S(G) $$
and a simple calculation shows that this composition is the identity $1_{S(G)}$. This proves claim 2.\\

Now consider the simplicial group $S(H)$ acting on the right on $S(G)$ and trivially
on $\Sigma$. The the map $\Omega'$ is $S(H)$-equivariant, hence we obtain a simplicial map
$\Omega=\Omega_a:\Sigma\times S(G)/S(H)\ra S(G)/S(H)$ and the following diagram of simplicial maps
commutes:
$$\begin{array}{ccc}
\Sigma\times S(G)&\stackrel{\Omega'}{\ra}&S(G)\\
\downarrow&&\downarrow\\
\Sigma\times S(G)/S(H)&\stackrel{\Omega}{\ra}&S(G)/S(H)
\end{array}$$

\noindent Write $\omega_a=|\Omega_a|$. Applying the geometric realization functor to the diagram above we
have the following commutative diagram:
$$\begin{array}{ccc}
\bS^1\times G^\bullet&\stackrel{\omega'_a}{\ra}&G^\bullet\\
\downarrow&&\downarrow\\
\bS^1\times G//H&\stackrel{\omega_a}{\ra}&G//H
\end{array}$$

\noindent Since $p:G^\bullet \ra G//H$ is onto and $p\alpha'$ represents $a$, using the diagram above
we conclude that $\omega_a$ satisfies the required properties. This proves Proposition 3.4.\vspace{.3in}

\noindent {\bf Proof of Proposition 3.2.} We just have to verify that the map 
$\pi_1(TOP(\bS^1\times\bS^{n-2})^\bullet)\ra\pi_1(TOP(\bS^1\times\bS^{n-2})//L)$
is onto. Since this map is a fibration it is enough to prove that $\pi_0(L^\bullet)\ra \pi_0(TOP(\bS^1\times\bS^{n-2})^\bullet)$
is one-to-one. Equivalently we have to prove that $\pi_0(L)\ra \pi_0(TOP(\bS^1\times\bS^{n-2}))$
is one-to-one. Note that $\pi_0(L)$ can be identified with $\pi_1(SO(n-1))$ by assigning to $[\alpha]\in \pi_1(SO(n-1))$
the component of $L$ containing $\hat{\alpha}\in L$ defined by $\hat{\alpha}(z,u)=(z,\alpha(z).u)$.
Let $p:\bS^1\times\bS^{n-2}\ra\bS^{n-2}$ denote projection onto the second factor. Recall that the Hopf construction associates
to each map $f:\bS^1\times\bS^{n-2}\ra\bS^{n-2}$ a map $H(f):\bS^n\ra\bS^{n-1}$ such that homotopic maps go to homotopic
maps (see \cite{Toda}, p.112). Also, the J-homomorphism \,\,$J:\pi_1(SO(n-1))\ra\pi_n(\bS^{n-1})$ is given by
$J([\alpha])=[H(p\,\circ\hat{\alpha})]$. But $J$ is one-to-one (see, for instance, \cite{Milnor-Kervaire}, p.512), and
\,$J$ factors through $\pi_0(TOP(\bS^1\times\bS^{n-2}))$ by the map $\pi_0(TOP(\bS^1\times\bS^{n-2}))\ra\pi_n(\bS^{n-1})$,
$[\phi]\mapsto [H(p\,\circ\phi)]$. Therefore the map $\pi_0(L)=\pi_1(SO(n-1))\ra \pi_0(TOP(\bS^1\times\bS^{n-2}))$
is also one-to-one. This proves Proposition 3.2.\\

Our proof of Proposition 3.3 depends on the following particular case of a result of Goodwillie \cite{Goodwillie} which is also
a consequence of Grunewald, Klein and Macko's Theorem 1.2 in \cite{GKM} (See also the
last Section of \cite{FO6}.)\\

\noindent {\bf Theorem.} {\it Let $p$ be an odd prime. Then $\pi_j A(\bS^1)$ is in $\cC_p$,
for $j<2p-2$.}\\

\noindent {\bf Proof of Proposition 3.3.}
As before write $N=\bS^1\times \bS^{n-2}$ and recall that in the proof of Theorem D we proved that
$\pi_j\Big( TOP(N)//L \Big)$ is $\cC_p$-isomorphic to $\pi_{j+1}\Big( \widetilde{TOP}(N)/TOP(N)  \Big)$,
for all $1\leq j$. Also recall that Hatcher's spectral sequence $E^r_{st}$
converges to $\pi_{s+t+1}\Big( \widetilde{TOP}(N)/TOP(N) \Big)$, for $s+t<k$ and that $E^2_{st}=H_s\Big(\Z_2;\pi_t\cP(N) \Big)$
is a subquotient of $\pi_t\cP(N)$. But $\pi_t\cP(N)\cong\pi_t\cP(\bS^1)$ (when $t<k$) and, from the Theorem above, $\pi_t\cP(\bS^1)$ is in
$\cC_p$, for all $t<k=2p-4$. Therefore $E^2_{st}$ is in $\cC_p$ for all $s,t$ with $s+t<k$. Consequently
$E^\infty_{st}$, $s+t<k$,  is in $\cC_p$ because these groups are subquotients of $E^2_{st}$.
Finally, since $\pi_{j+1}\Big( \widetilde{TOP}(N)/TOP(N) \Big)$ has a finite length filtration with successive quotient
groups $E^\infty_{st}$, $j=s+t$ and $\cC_p$ is a Serre class, it follows that 
$\pi_{j+1}\Big( \widetilde{TOP}(N)/TOP(N) \Big)$ is in $\cC_p$, for $j<k$. This proves Proposition 3.3.\\
\vspace{1in}

\noindent {\bf \Large  Section 4. Proof of Proposition 3.1.}\\

We will use the following facts about the ring $\zp$.\\

\begin{enumerate}
\item[{\it (i)}] If $0\ra A\ra B\ra C\ra  0$ is an exact sequence of abelian groups, then 
$0\ra A\ozp\ra B\ozp\ra C\ozp\ra  0$ is an exact sequence of $\zp$-modules.
(This is because $\zp$ is torsion free.)

\item[{\it (ii)}] A submodule of a finitely generated $\zp$-module is finitely generated as a $\zp$-module.
(This is because $\zp$ is a principal ideal domain.) 

\item[{\it (iii)}] For a prime $q\neq p$ and any abelian group $A$, we have $\tau_q(A)\ozp =0$.

\item[{\it (iv)}] For any abelian group $A$ we have that $\tau_p(A)\ozp =\tau_p(A)$.
(Proof: Clearly $\Z_{p^n} \ozp\cong\Z_{p^n}$. Applying this and (i) to the subgroup generated by a supposed element in the kernel of
$\tau_p(A)\ra\tau_p(A)\ozp $, we get
that $\tau_p(A)\ra\tau_p(A)\ozp $ is monic. Any element in
$\tau_p(A)\ozp $ can be written in the form $\frac{a}{s}=a\otimes \frac{1}{s}$, $(s,p)=1$, $a\in A$, $p^na=0$, for some $n$. 
Hence there are integers $\lambda $ and $\mu$ such that $\lambda p^n+\mu s=1$. Then $\frac{a}{s}=\frac{(\lambda p^n+\mu s)a}{s}=\mu a$.
Therefore $\frac{a}{s}$ is in the image of $\tau_p(A)\ra\tau_p(A)\ozp $.)

\item[{\it (v)}] If $C$ is a finitely generated $\zp$-module, then $C$ is isomorphic to a finite sum, where each summand is either
$\zp$ or $\zp/p^n\zp=\Z/p^n\Z=\Z_{p^n}$, for some $n$. This is because $\zp$ is a principal ideal domain. 
\end{enumerate}

Recall that an abelian group $A$ is in the class $\cC_p$ if
$A\otimes \zp$ is finitely generated as a $\zp$-module.
If $A$ is in $\cC_p$ then, by items (i), (ii), (iv) and (v) $\tau_p(A)\ozp =\tau_p(A)$ is finitely generated, hence finite.
Therefore, since by (i) $ 0\ra \tau(A)\otimes \zp\ra A\otimes \zp\ra \Big( A/\tau(A)\Big)\otimes \zp\ra 0$
is exact, we have that 
$A$ is in $\cC$ if and only if $\tau_p(A)$ is finite and $\Big( A/\tau(A)\Big) \otimes \zp$ is finitely generated 
as a $\zp$ module. Hence, for $A$ a torsion group, $A$ being in $\cC_p$ is equivalent to $\tau_p(A)$ being finite.
On the other hand, for $A$ torsion free, $A$ being in $\cC_p$ is equivalent to $A$ being (isomorphic to) a subgroup of $(\zp)^k=\zp\oplus...\oplus\zp$,
for some $k$. This follows from (v) and the fact that $A\ra A\ozp$ is injective. (The map $A\ra A\ozp$ is injective
because $\Z\ra\zp$ is injective and $A$ is torsion free.)\\


We have to prove (see Spanier \cite{Sp}, chap. 9, sec. 6):
\begin{enumerate}
\item[{\bf (a)}] {\it $\cC_p$ contains the trivial group.}

\item[{\bf (b)}] {\it If $ A$ is in $\cC_p$ and $A'$ is isomorphic to $A$, then $A'$ is in $\cC_p$.}

\item[{\bf (c)}] {\it If $A$ is in $\cC_p$ and $B\sbs A$, then   $B$ is in $\cC_p$.}

\item[{\bf (d)}] {\it If $A$ is in $\cC_p$ and $B\sbs A$, then   $A/B$ is in $\cC_p$.}

\item[{\bf (e)}] {\it If $\,\,\, 0\ra A\ra B\ra C\ra 0$ is a short exact sequence and $A,C$ are in $\cC_p$, then $B$ is in $\cC_p$.}

\item[{\bf (f)}] {\it If $\,\,A, B$ are in $\cC_p$, then $A\otimes B$ is in $\cC_p$.}

\item[{\bf (g)}] {\it If If $\,\,A, B$ are in $\cC_p$, then $Tor(A,B)$ is in $\cC_p$.}

\item[{\bf (h)}] {\it If $\,\,A$ is in $\cC_p$, then $H_j(A)$ is in $\cC_p$, $\, j>0$.}
\end{enumerate}


Properties (a) and (b) are obviously true. Property (c) follows from facts (i) and (ii) above.\\

We prove (d). Let $A$ be in  $\cC_p$ and $B\sbs A$.
Since $0\ra B\ra A\ra (A/B)\ra 0$ is exact, by fact (i) we obtain that
$A\ozp\ra\Big( A/B\Big)\ozp$ is onto. Hence $\Big( A/B\Big)\ozp$ is finitely generated. This proves (d).\\

Property (e) follows directly from fact (i).\\

We prove (f). Let $A$ and $B$ be in $\cC_p$. Since $\zp\ozp\cong \zp$ we have that
$(A\otimes B)\ozp$ and $(A\ozp)\otimes (B\ozp)$ are isomorphic as $\zp$-modules.
It follows that $(A\otimes B)\ozp\cong (A\ozp)\otimes (B\ozp)$
is finitely generated as a $\zp$-module.\\

To prove (g) let $A$ and $B$ be in
$\cC_p$ and
recall that an abelian group $C$ is in $\cC_p$ if and only if $\tau_p(C)$ is finite and
$\Big( C/\tau(C) \Big)\ozp$ is finitely generated as a $\zp$-module. Therefore, since $Tor(A,B) $ is a torsion group,
to prove (g) we just have to prove that $\tau_p\Big( Tor(A,B) \Big)$ is finite.
But $Tor(A,B)=Tor\Big( \tau_p(A),\tau_p(B)\Big)\oplus\bigoplus_{q\neq p}Tor\Big( \tau_q(A),\tau_q(B)\Big)$.
Therefore $\tau_p\Big( Tor(A,B) \Big)=Tor\Big( \tau_p(A),\tau_p(B)\Big)$, which is finite
because  $\tau_p(A)$ and $\tau_p(B)$ are finite.\\

Finally, we prove (h). First note that, by taking the homology long exact sequence induced by
$0\ra \tau(A)\ra A\ra\Big(A/\tau(A) \Big)\ra 0$, we can see that proving  (h) is equivalent to proving (h) in
the following two special cases: when $A$ is a torsion group, and when $A$ is torsion free.\\

Let $A$ be a torsion group in $\cC_p$. Then $A=\tau_p(A)\oplus\Bigg( \bigoplus_{q\neq p}\tau_q(A)\Bigg)$.
Write $B=\bigoplus_{q\neq p}\tau_q(A)$. Since $\tau_p(A)$ is a finite $p$-group, it is in the acyclic Serre
rings 5 and 7 of p. 505 of Spanier \cite{Sp}. Hence $H_j\Big(\tau_p(A)\Big)$ is also a finite $p$-group,
for all $j>0$. Also, since $B$ is a torsion group with no $p$-torsion, it is in the acyclic Serre ring
8 of p. 505 of Spanier. Therefore $H_j\Big(B)$ is also a torsion group with no $p$-torsion. These facts,
together with the K\"unneth Formula imply that $H_j(A)=H_j\Big(\tau_p(A)\Big)\oplus H_j(B)$, for all $j>0$.
Then, by facts (iii) and (iv) above, $H_j(A)\ozp=H_j\Big(\tau_p(A)\Big)$, which is finite. This proves (h) when $A$ is a torsion group.\\

Let $A$ be a torsion free group in $\cC_p$. Then $A$ is a subgroup of $(\zp)^k$, for some $k$.
Note that any finitely generated subgroup of $(\zp)^k$ is a free group of rank at most $k$. 
Ordering properly the products of powers of primes $q$, $q\neq p$, we can find a sequence of integers
$s_0, s_1,s_2 ...$    such that $s_i| s_{i+1}$, $p\not\,\mid s_i$, and 
in addition satisfying the following condition: for every $s$, with $p\not\,\mid s$, there is $i$
such that $s|s_i$. Define the free rank k subgroup $B_i$ of $(\zp)^k$ as $B_i=\frac{1}{s_i}\Z^k$.
Then $(\zp)^k=\bigcup B_i$. Therefore $(\zp)^k=\lim_{\ra}B_i$. Note that $B_{i+1}/B_i=(\Z_{m})^k$, with $m=\frac{s_{i+1}}{s_i}$,
thus the finite groups $B_{i+1}/B_i$ have no $p$-torsion. Hence, we also have that
$A=\lim_{\ra}\Big( A\cap B_i\Big)$, and $\Big( A\cap B_{i+1}\Big)\, /\,\Big( A\cap B_i\Big)$
has no $p$-torsion. Since homology and tensor products commute with direct limits we have that
$H_j(A)\ozp=\lim_{\ra} H_j\Big( A\cap B_i\Big)\ozp$.
We claim the maps of the direct system
$ \Big\{ H_j\Big( A\cap B_i\Big)\ozp\Big\}$ are all isomorphisms. Consequently  $H_j(A)\ozp\cong H_j(\Z^\ell)\ozp$, for some $\ell$
and all $j>0$.  And this proves (h) once we verify our claim. For this note that $H_j(A\cap B_i)\ozp=H_j(A\cap B_i,\,\zp)$ by the Universal
Coefficient Theorem. Next apply Proposition 9.5 (ii) of \cite{Brown}, p. 82, in which we specify $H=A\cap B_i$ and $G=A\cap B_{i+1}$
to conclude that $\varphi :H_j(A\cap B_i,\,\zp)\ra H_j(A\cap B_{i+1,}\,\zp)$ is onto because the index $[G:H]$ is invertible in
$\zp$. Consequently $\varphi$ is an isomorphism since $\zp$ is a principal ideal domain and the domain and range of $\varphi$ are
isomorphic finitely generated free $\zp$-modules.
\vspace{.7in}

\noindent {\bf \Large  Section 5. Proof of Lemma 2.3.1.}\\

The proof of Lemma 2.3.1 depends on the following facts.\\

\noindent {\bf Fact 1.} For each finite subgroup $G$ of $\pi_k\Big(P(\bS^1\times\bS^{n-2})\Big)$,
there exists a positive integer $j$ such that
$$(\bar\nu_j)_* (G)=0\,\,\,\,\,\,\,\mbox{in} \,\,\,\,\pi_k\Big(P(\bS^1\times\bS^{n-2})\Big)$$
where $\bar\nu_j: P(\bS^1\times \bS^{n-2})\ra P(\bS^1\times \bS^{n-2})$ is induced by the $j$-sheeted cover 
$\bS^1\times\bS^{n-2}\ra \bS^1\times\bS^{n-2}$ in the same way $\nu_j$ is induced.\\

\noindent {\bf Fact 2.} For each positive integer $j$ there is a natural commutative diagram:

$$\begin{array}{ccc}
\pi_k\Big(P^s(\bS^1\times\bS^{n-2})\Big)&\stackrel{F_*}{\ra}&\pi_k\Big(P(\bS^1\times\bS^{n-2})\Big)\\
(\nu_j)_*\uparrow&&\uparrow(\bar\nu_j)_*\\
\pi_k\Big(P^s(\bS^1\times\bS^{n-2})\Big)&\stackrel{F_*}{\ra}&\pi_k\Big(P(\bS^1\times\bS^{n-2})\Big)
\end{array}$$\\
where $F:P^s (\bS^1\times\bS^{n-2}) \ra P(\bS^1\times\bS^{n-2})$ is the forget structure map.\\

\noindent {\bf Fact 3.}  $ker\, F_*$ is a finitely generated abelian group provided 
$k<<n$.\\

Before justifying these facts, we use them to prove Lemma 2.3.1.\\

Because of Fact 3, the torsion subgroup $T_k$ of $ker \, F_*$ is finite. Let $S$ be a finite subgroup
of $T$ with $|S|>|T_k|$. Note that such $S$ exist with arbitrarily large finite cardinality 
since $T$ is a torsion abelian group of infinite cardinality. (T is obviously abelian when $k \geq 1.$
The case $k=0$ follows from the fact that the stable smooth pseudo-isotopy space 
${\cP}^s (\bS^1\times\bS^{n-2})$ is an infinite loop space; see \cite{Hat}, Appendix II.)\\

Because of Facts 1 and 2, there exists a positive integer $j$ such that 
$$(\nu_j)_* (S)\subseteq T_k.$$

Consequently the cardinality of the subgroup $ker \Big( (\nu_j)_*|_S\Big)$ of $\overline T$
is at least as large as $|S|/|T_k|$. This completes the proof of Lemma 2.3.1 (modulo verifying Facts 1-3)
 since $|S|/|T_k|$ can be arbitrarily large.\\

\noindent {\bf Verification of Facts 1-3.}

\noindent Fact 2 needs no justification since it is obvious.\\
Fact 3 is a consequence of results of Burghelea-Lashof \cite{BL}, see
 \cite{Hat}, Theorem 5.5 and Corollary 5.6, in conjunction with Dywer's result \cite{Dw} that
 $\pi_i\Big( P^s(\D^n)\Big)$ is finitely generated for $i<<n$.\\
 Fact 1 can be deduced from Quinn's paper \cite{QuinnII} by observing that $\nu_j(f)$ becomes as 
 controlled over $\bS^1$ (as needed) as $j \ra \infty$; then use the result that
 $\pi_i\Big( P(\bS^{n-2})\Big)=0$ for $i<<n$ which is a consequence of \cite{Hat}, Corollary 5.5, together with the
 well known fact that $P(\D^n)$ is contractible (via the Alexander isotopy).
 \vspace{.8in}

\noindent {\bf \Large  Section 6. Proof of Proposition  2.1.5.}\\

Consider the following diagram

$$\begin{array} {ccccccccccccc}
S&&\stackrel{\alpha}{\longrightarrow}&&K&&\stackrel{\Omega_\gamma}{\longrightarrow}&&
\met_\gamma^\epsilon(M)&&\stackrel{\Lambda^\epsilon_\gamma}{\longrightarrow}&&
P(\bS^1\times\bS^{n-2})\\
\\
&&&& &&{\scriptstyle{\Omega'_\gamma}}\searrow &&\downarrow {\mbox{lift}} && \,\,\,\,\,\,\,\,\,\,\,\nearrow{\scriptstyle{\Lambda^\epsilon_\gamma}}\\
 \\
&&&&&&&&\met_\gamma^\epsilon(Q)

\end{array}$$\\

\noindent where we wrote a ``prime" on $\Omega'_\gamma:K\ra\met_\gamma^\epsilon(Q)$ to differentiate it from
$\Omega_\gamma:K\ra\met_\gamma^\epsilon(M)$. Observe that the right hand triangle in the diagram above is commutative 
(see Remark before Lemma 2.1.3).\\

\noindent {\bf Claim 1.}
{\it $\Lambda_\gamma^\epsilon\Omega_\gamma'$ is homotopic to the inclusion $\iota:K\ra P(\bS^1\times\bS^{n-2})$.}\\

\noindent{\bf Proof of Claim 1.} This is a consequence of the fact that
$\Lambda_\gamma^\epsilon\Omega_\gamma'(\varphi)$ is equal to $\varphi$, up to rescaling in the $t$-direction.
And this follows from the definitions of  $\Lambda_\gamma^\epsilon$ and $\Omega_\gamma'$ and the following
fact: ``let $g_t$ be a family of Riemannian metric on $M$ ($g_t(x)$ smooth with respect to $(x,t)$) and define the Riemannian metric
$\bar{g}(x,t)=g_t(x)+dt^2$ on $M\times \R$. Then the vertical lines $t\mapsto (x,t)$ are geodesics". This is a consequence of
Koszul's formula. This proves Claim 1.\\

\noindent {\bf Claim 2.} (lift) $\circ\,\Omega_\gamma\circ \alpha$\, \, {\it is homotopic to \,$\Omega_\gamma'\circ \alpha$.}\\

\noindent {\bf Proof of Claim 2.} We prove this by giving two deformations.\\

\noindent {\bf First deformation.}
By hypothesis, there is a deformation $h_v$, $v\in [0,1]$, with $h_0(u)=\nu_j(\alpha(u))$, $u\in S$, and 
$h_1(u)=1_{\bS^1\times\bS^{n-2}\times [0,1]}$.
Let $K'=\{ h_v(u)\,\, |\,\, u\in S,\,\, v\in [0,1]\}$. Then $K'$ is compact and let $r'>0$ be given by 
Lemma 1.2.2 for $K'$ and $j\ell$. We will assume that $r>r'$.\\

Let $R$ denote the normal tubular neighborhood of $\gamma\sbs M$ of width 6$r$.
The lift of $\gamma\sbs M$ to $\bS^1\times\R^{n-1}$ gives, as mentioned in Section 2.1,
exactly one closed geodesic, which we call also $\gamma\sbs \bS^1\times\R^{n-1}$; but also gives
a countable number of disjoint infinite geodesic lines $\ell_1, \ell_2,...$. In the same way lifting
the neighborhood $R$ we obtain one copy of $R$ (which contains $\gamma$ and we also denote by $R$)
plus disjoint neighborhoods $R_i$'s of the $\ell_i$'s each diffeomorphic to a cylinder $\R\times\D^{n-1}$.\\

Let $u\in S$ and write $\omega_u=\Omega_\gamma (\alpha(u))$. 
Denote by $\bar{\omega}_u$ the Riemannian metric which is the lift of $\omega_u$ to $Q=\bS^1\times\R^{n-1}$,
i.e. $\bar{\omega}_u={\mbox{lift}}(\omega_u)$. Note that
$\bar{\omega}_u$ outside $R\cup\bigcup R_i$ coincides with the hyperbolic metric $\rho$.
Also $\bar{\omega}_u$ on $R$ coincides with $\omega_u$, and
on each $R_i$ it is the lifting of $\omega_u$ on $R\sbs M$. Since the cover $R_i\ra R$ 
can be written as a composition $R_i\stackrel{p_i}{\ra} R \stackrel{p}{\ra}R$, where
$R\stackrel{p}{\ra} R$ the $j$-sheeted cover,  we have that $\bar{\omega}_u$ on each $R_i$ is also the lifting, by $p_i$, of the metric
$\Omega''_\gamma (\nu_j(\alpha(u)))$, where we put the ``two primes" on $\Omega''_\gamma:\nu_j(K)\ra\met_\gamma^\epsilon (R)$,
to differentiate it from the other $\Omega_\gamma$'s.
(Here $R$ is considered with the pulled back metric by the $j$-sheeted cover $R\ra R$. Note that, even though $R$ is not complete, it still
makes perfect sense to define $\Omega'_\gamma$ as in section 1.2, provided $r$ is large enough.)\\

Finally, define $\Big(\bar{\omega}_u\Big)_v$ as being equal to $\bar{\omega}_u$ outside $\bigcup R_i$
and equal to the lifting $(p_i)^*\Omega_\gamma'' (h_v(u))$ on each $R_i$. In this way we can deform
$\bar{\omega}_u$ to $\hat{\omega}_u$, which is equal to $\rho$ outside $R\cup\bigcup R_i$, equal to $\omega_u$ on $R$ and
equal to $\rho''$ on each $R_i$, after properly identifying each $R_i$ with $\R\times D^{n-1}$, where $D^{n-1}$ is a disc
of large enough radius. (To recall the definition of $\rho''$ see Section 1.2.)
This completes the construction of the first deformation.\\

\noindent {\bf Second deformation.}\\ 
We deform now $\hat{\omega}_u$ to $\omega'_u=\Omega'_\gamma(\alpha(u))$. To do this just define
$\Big(\hat{\omega}_u\Big)_v$, $v\in [0,1]$, to be equal to $\hat{\omega}_u$ outside $\bigcup R_i$ and on each $R_i$ equal to
$\rho''_v$, where $\rho''_v$ is the deformation given in the remark after Lemma 1.2.1.\\

To prove that these deformations are continuous and their images lie in $\met (Q)$ (recall that $\met(Q)$ is really
$\met(Q,\rho_0)$ with the smooth quasi-isometry topology, see section 1.4) just note the following two facts:\\

\noindent {\bf (i).} The first deformation only happens on each $R_i$, and all metrics
$\Big(\bar{\omega}_u\Big)_v$ together with $\rho_0$ are invariant by a cocompact action (of translations in the $\ell_i$
direction of certain length) coming from the cover $p_i:R_i\ra R$. Now just apply Lemma 1.4.1.\\

\noindent {\bf (ii).} For the second deformation we proceed as in the previous case. Just note that $\rho''_v$ is actually invariant
by any translation in the infinite line direction. \\

This completes the construction of the second deformation and  the proofs of Claim 2 and Proposition 2.1.5.

\vspace{.8in}

\noindent {\bf \Large  Appendix}\\

The following Lemma is needed in the last part of the proof of Proposition 1.1.1.\\

\noindent {\bf Lemma.} {\it Let $P$ be a closed smooth $k$-manifold and $F:P\times\R^m\ra P\times\R^m$ a smooth embedding
such that $F(p,0)=(p,0)$ and $DF(p,0)$ is the identity, for all $p\in P$. Then $F$ is smoothly isotopic to the identity $1_{P\times\R^m}$,
relative to $P\times \{ 0\}$.}\\

\noindent {\bf Proof.} Write $F=(f,g)$, $f:P\times\R^m\ra P$, $g:P\times\R^m\ra\R^m$. 
Since $g(p,0)=0$ there are smooth maps $h_i:P\times\R^m\ra \R$ such that:
\begin{enumerate}
\item[{\bf 1.}] $g(p,q)=q_1h_1(p,q)+...+q_mh_m(p,q)$, where $q=(q_1,...,q_m)\in\R^m$.

\item[{\bf 2.}] $h_i(p,0)=e_i$, where the $e_i$'s form the canonical basis of $\R^m$.
\end{enumerate}

Define $$F_t(p,q)=\Big( f(p,tq),\frac{1}{t}\, g(p,tq)\Big) = \Bigg( f(p,tq),\sum_i q_ih_i(p,tq)\Bigg)$$

\noindent where the last equality is given by item {\bf 1} above.
Then $F_t$ is smooth for all $(p,q,t)$ and it is clearly an embedding, for all $t\neq 0$. For $t=0$, by item {\bf 2}
we have $F_0(p,q)=\Big ( f(p,0), \sum_iq_ih_i(p,0)\Big)=(p,q)$. Since $F_1=F$, this completes the proof of the Lemma.

F.T. Farrell

SUNY, Binghamton, N.Y., 13902, U.S.A.\\

P. Ontaneda

SUNY, Binghamton, N.Y., 13902, U.S.A.

\end{document}